\documentclass[review]{elsarticle}
\usepackage{amsmath}
\usepackage{amssymb}
\usepackage{amsthm}
\usepackage{caption}
\usepackage{lineno,hyperref}
\usepackage{subfigure}
\usepackage{graphics}
\usepackage{booktabs}
\usepackage{color}
\modulolinenumbers[5]

\newtheorem{lemma}{Lemma}

\newtheorem{thm}{Theorem}
\newtheorem{cor}{Corollary}
\newtheorem{remark}{Remark}

\begin{document}

\begin{frontmatter}

\title{A necessary and sufficient condition for lower bounds on crossing numbers of generalized periodic graphs in an arbitrary surface}

%
%
%
\author[label1]{Xiwu Yang\corref{cor1}}
\ead{cinema@lnnu.edu.cn}
\author[label1]{Xiaodong Chen}
\author[label2]{Yuansheng Yang}

\address[label1]{School of Mathematics, Liaoning Normal University, Dalian 116029,  P.R. China}
\address[label2]{Department of Computer Science, Dalian University of Technology, Dalian 116024, P.R.China}
\cortext[cor1]{Corresponding author}

\begin{abstract}
Let $H$, $T$ and $C_n$ be a graph, a tree and a cycle of order $n$, respectively.
Let $H^{(i)}$ be the complete join of $H$ and an empty graph on $i$ vertices.
Then the Cartesian product $H\Box T$ of $H$ and $T$ can be obtained by applying zip product on $H^{(i)}$ and the graph produced by zip product repeatedly.
Let $\textrm{cr}_{\Sigma}(H)$ denote the crossing number of $H$ in an arbitrary surface $\Sigma$.
If $H$ satisfies certain connectivity condition,
then $\textrm{cr}_{\Sigma}(H\Box T)$ is not less than the sum of the crossing numbers of its ``subgraphs".
In this paper, we introduced a new concept of generalized periodic graphs, which contains $H\Box C_n$.
For a generalized periodic graph $G$ and a function $f(t)$,
where $t$ is the number of subgraphs in a decomposition of $G$,
we gave a necessary and sufficient condition for $\textrm{cr}_{\Sigma}(G)\geq f(t)$.
As an application, we confirmed a conjecture of Lin et al. on the crossing number of the generalized Petersen graph $P(4h+2,2h)$ in the plane.
Based on the condition, algorithms are constructed to compute lower bounds on the crossing number of generalized periodic graphs in $\Sigma$.
In special cases, it is possible to determine lower bounds on an infinite family of generalized periodic graphs, 
by determining a lower bound on the crossing number of a finite generalized periodic graph.
\end{abstract}
\end{frontmatter}


\section{Introduction}
In this paper, we only consider simple graphs if not specified.
For some terminology and notation not defined here,
readers can refer to \cite{Bon}.
A drawing of a graph is \emph{good}, provided that no edges crosses itself,
no adjacent edges cross each other, no two edges cross more than once,
and no three edges cross at a common point.
In a good drawing, a common point of two edges other than endpoints is a \emph{crossing}.
The crossing number of a graph $G$ in an arbitrary surface $\Sigma$, denoted by $\textrm{cr}_{\Sigma}(G)$, is the minimum number of
crossings among all good drawings of $G$ in $\Sigma$.
In particular, let $\textrm{cr}(G)$ denote the crossing number of a graph $G$ in the plane.

When determining the crossing number of a graph, it is natural to decompose the graph to obtain simpler instances.
Let $G=(V,E)$ be a connected graph and let $F\subseteq E(G)$ be a minimum edge cut in $G$ of size $d=|F|$.
Let $H_1$ and $H_2$ be the two components of $G-F$, and we define $G_i$
obtained from $H_i$ by adding a new vertex $v_i$ and connecting it to all the endvertices of $F$ in $H_i$, $1\leq i\leq 2$.
It is shown that $\textrm{cr}_{\Sigma}(G)=\textrm{cr}_{\Sigma}(G_1)+\textrm{cr}_{\Sigma}(G_2)$ if $1\leq d\leq 3$ \cite{Bok3,Lea}.
If $d\geq 4$, there are graphs that $\textrm{cr}(G)<\textrm{cr}(G_1)+\textrm{cr}(G_2)$ \cite{Bea}.

Assume that $F$ is a matching of $G$.
Let $N_i=N_{G_i}(v_i)$ be the set of neighboring vertices of $v_i$ in $G_i$, $1\leq i\leq 2$,
and let $\sigma: N_1\rightarrow N_2$ be a bijection.
We call $\sigma$ a \emph{zip function} of the graphs $G_1$ and $G_2$ at vertices $v_1$ and $v_2$.
The \emph{zip product} of $G_1$ and $G_2$ according to $\sigma$ is the graph $G_1\odot_{\sigma}G_2$
obtained from the disjoint union of $G_1-v_1$ and $G_2-v_2$ after adding edges $u\sigma(u)$ for any $u\in N_1$.
Then $G$ is a zip product of $G_1$ and $G_2$.
For any bijection $\sigma: N_1\rightarrow N_2$, it is proved that if both $H_1$ and $H_2$ satisfy certain connectivity condition \cite{Bok,Bok3},
then $\textrm{cr}_{\Sigma}(G_1\odot_{\sigma} G_2)\geq \textrm{cr}_{\Sigma}(G_1)+\textrm{cr}_{\Sigma}(G_2)$;
if both $H_1$ and $H_2$ satisfy other connectivity conditions \cite{Bok07,Ouy},
then $\textrm{cr}(G_1\odot_{\sigma} G_2)\geq \textrm{cr}(G_1)+\textrm{cr}(G_2)$.

Let $H$ be a graph and let $T$ be a tree.
Let $H^{(j)}$ be the complete join of $H$ and an empty graph on $j$ vertices. 
Then the Cartesian product $H\Box T$ of $H$ and $T$ can be obtained by applying zip product on $H^{(j)}$ and the graph produced by zip product repeatedly.
If $H$ satisfies certain connectivity condition,
then $\textrm{cr}(H\Box T)$ is not less than the sum of the crossing numbers of its ``subgraphs" \cite{Bok,Bok07,Ouy}.
This lower bound on $\textrm{cr}(H\Box T)$ confirmed \cite{Bok} the conjecture of Jendrol' and \v{S}\v{c}cerbov\'{a} \cite{Jen} 
and determined \cite{Bok07,Ouy} crossing numbers of several graph families.

Computer programs are also applied to determine the crossing number of a graph in the plane.
Lin et al. \cite{Lin09} obtained an algorithm called CCN which can calculate the crossing number of any graph of small order.
In \cite{Lin09}, they conjectured that $\textrm{cr}(P(4h+2,2h))=2h+1$, $h\geq 3$, where $P(4h+2,2h)$ is a generalized Petersen graph.
Chimani et al. \cite{Chi} gave an algorithm which can compute the crossing number of any graph with small order.
For any maximal pathwidth-3 graph $G$, Biedl et al. \cite{Bie} proved that $\textrm{cr}(G)$ equals the rectilinear crossing number of $G$,
and proposed a linear time algorithm to compute $\textrm{cr}(G)$.

Let $C_n$ be a cycle of order $n$.
By applying zip product on $H^{(2)}$ and the graphs produced by zip product repeatedly, we can obtain $H\Box C_n$.
Can we obtain a lower bound on $\textrm{cr}_{\Sigma}(H\Box C_n)$ similarly?
In this paper, we answer this question affirmatively.
We introduced a new concept of a generalized periodic graph, which is a natural generalization of a periodic graph \cite{RefDvo,RefPin}.
For a generalized periodic graph $G$ and a function $f(t)$,
where $t$ is the number of subgraphs in a special decomposition of $G$,
we gave a necessary and sufficient condition for $\textrm{cr}_{\Sigma}(G)\geq f(t)$.
As an application, we confirmed the conjecture from Lin et al. \cite{Lin09}.
Our condition can be applied to construct algorithm to compute the lower bounds on the crossing number of generalized periodic graphs in $\Sigma$.
For example, we give algorithms to compute the lower bounds on the crossing number of generalized periodic graphs in the plane.
In special cases, it is possible to determine lower bounds on an infinite family of generalized periodic graphs, 
by determining a lower bound on the crossing number of a finite generalized periodic graph.

\section{Generalized periodic graphs}
Loops and parallel edges are allowed in the graphs of this section.
A periodic graph was introduced by Pinontoan and Richter \cite{RefPin}.
In this paper, we use the definition proposed by Dvo\v{r}\'{a}k and Mohar \cite{RefDvo}.
A \emph{tile} is a triple $Q=(G,L,R),$ where $G$ is a graph, $L=(l_1,...,l_k)$ and
$R=(r_1,...,r_k)$ are sequences of vertices in $G$ of the same length $k$.
Each vertex in $L$ or $R$ can repeatedly appear.
The length of each sequence is called the \emph{width} of the tile.
Suppose that $G_1$ and $G_2$ are two vertex disjoint graphs
and $Q_1=(G_1,L_1,R_1)$ and $Q_2=(G_2,L_2,R_2)$ are tiles of the same width $k$.
Let $Q_1Q_2$ denote the tile $(H,L_1,R_2)$, where $H$ is the graph obtained from the vertex disjoint union of $G_1$ and $G_2$ by adding
an edge between the $j$-th vertices of $R_1$ and $L_2$, $1\leq j\leq k$.
Assume that $Q_1^1=Q_1$ and $Q_1^t=Q_1^{t-1}Q_1$ for integer $t>1$.
Let $Q_1^t=(G_t,L_1,R_t),t\geq 1$.
A \emph{periodic graph} $\odot(Q_1^t),t\geq 2$, is a graph obtained from $G_t$ by adding edges between the $j$-th vertices of $L_1$ and $R_t$ for $j=1,...,k$.
The edges of $\odot(Q_1^t)$ that belong to the copies of $G_1$ are \emph{internal}, and the edges between the copies are \emph{external}.
It is evident that the order of a periodic graph is not a prime.
Let $\odot(Q^n)$ be a periodic graph.
Pinontoan and Richter \cite{RefPin} proved that
$$c(Q)={\lim_{n \to +\infty}}\textrm{cr}(\odot(Q^n))/n$$ exists.
Dvo\v{r}\'{a}k and Mohar \cite{RefDvo} showed that for every $\varepsilon>0$,
there exists a computable constant $N=O(1/\varepsilon^6)$ such that $|\textrm{cr}(\odot(Q^t))/t-c(Q)|\leq \varepsilon$ for every integer $t\geq N$.

Let $G$ be a simple graph.
A \emph{decomposition} of $G$ is a list of edge-disjoint subgraphs of $G$ such that each edge of $G$ appears in one of the subgraphs in the list.
A decomposition $\{H_1,...,H_t\}$ of $G$ is \emph{transitive}
if for each pair of integers $i$ and $j, 1\leq i,j\leq t$,
there is an automorphism $\theta_{i,j}$ of $G$ such that
$uv\in E(H_{i+\ell})$ if and only if $\theta_{i,j}(u)\theta_{i,j}(v)\in E(H_{j+\ell}), 0\leq \ell< t$, where the subscripts are taken modulo $t$.
A graph $G$ is \emph{generalized periodic} if there exists a transitive decomposition of $G$.

In a periodic graph $\odot(Q_1^t)$ $(t\geq 3)$, let $S_i$ be the subgraph of $\odot(Q_1^t)$ induced by the set of all external edges between
the $i$-th and the $(i+1)$-th copy of $G_1$ for $1\leq i\leq t-1$,
and let $S_t$ be the subgraph of $\odot(Q_1^t)$ induced by the set of all external edges between the first and the $t$-th copy of $G_1$.
Let $G_i^+$ be the union of $S_i$ and the $i$-th copy of $G_1$ for $1\leq i\leq t$.
It is easy to check that $\{G_1^+,...,G_t^+\}$ is a transitive decomposition of $\odot(Q_1^t)$.
Then each periodic graph $\odot(Q_1^t)$ $(t\geq 3)$ is generalized periodic.

A generalized periodic graph may be not periodic.
Let $V(K_n)=\{v_i:1\leq i\leq n\}, n\geq 3$.
Suppose that $n=2r+1$.
Let $V(H_i)=\{v_{i+j}:0\leq j \leq r\}$ and $E(H_i)=\{v_iv_{i+j}:1\leq j\leq r\}, 1\leq i\leq 2r+1$, where the subscripts are taken modulo $2r+1$.
It is easy to check that $\{H_1,...,H_{2r+1}\}$ is a transitive decomposition of $K_{2r+1}$,
and then $K_{2r+1}$ is generalized periodic.
If $2r+1$ is a prime,
then $K_{2r+1}$ is not a periodic graph.

\section{Lower bounds on crossing numbers of generalized periodic graphs in an arbitrary surface}
Let $G$ be a graph and let $A$ and $B$ be subgraphs of $G$.
In a good drawing $D$ of $G$ in an arbitrary surface $\Sigma$, the number of crossings crossed by one edge in $A$
and the other edge in $B$ is denoted by $\textrm{cr}_D(A,B)$.
Especially, $\textrm{cr}_D(A,A)$ is denoted by $\textrm{cr}_D(A)$ for short.
The following two equalities are trivial.
\begin{eqnarray}
\textrm{cr}_D(A\cup B)&=&\textrm{cr}_D(A)+\textrm{cr}_D(B)+\textrm{cr}_D(A,B),\label{I1}\\
\textrm{cr}_D(A,B\cup C)&=&\textrm{cr}_D(A,B)+\textrm{cr}_D(A,C),\label{I2}
\end{eqnarray}
where $A$, $B$ and $C$ are pairwise edge-disjoint subgraphs of $G$.
The number of crossings of $D$, denoted by $\textrm{cr}(D)$, is $\textrm{cr}_D(G)$.
For a graph $G$, let $H$ be a subgraph of $G$ and let $f_D(H)$ be a mapping from $H$ to the set of all nonnegative real numbers:
\begin{equation}\label{fun}
f_D(H)=\textrm{cr}_D(H)+\textrm{cr}_D(H,G\setminus E(H))/2.
\end{equation}

Let $\{H_1,...,H_t\}$ be a transitive decomposition of a graph $G$.
Assume that $H_{i,j}$, $j-i<t$, is the union of $H_i,...,H_j$ with the subscripts taken modulo $t$.
Suppose that $D$ is a good drawing of $G$.
Let $c$ be a positive real number.
For $1\leq i\leq t$,
if there exists a positive integer $l_i$ satisfying
$f_D(H_{i,i-1+l_i})\geq l_ic$,
then let $l^D_c(H_i)$ be the smallest $l_i$ such that $f_D(H_{i,i-1+l_i})\geq l_ic$,
and $l^D_c(H_i)=t+1$, otherwise.

\begin{lemma}\label{decom}
Suppose that a graph $G$ is generalized periodic with a transitive decomposition $\{H_1,...,H_t\}$
and $D$ is a good drawing of $G$ in an arbitrary surface $\Sigma$.
Let $c$ be a positive real number.
If max$\{l^D_c(H_1),...,l^D_c(H_t)\}\leq t$,
then there exists a decomposition $\{G_1,...,G_k\}$ of $G$,
where $G_i=H_{r_i,r_i-1+h_i}$
and $h_i$ is the smallest positive integer such that
$f_D(H_{r_i,r_i-1+h_i})\geq ch_i,1\leq i\leq k$.
\end{lemma}
\begin{proof}
Without loss of generality,
assume that max$\{l^D_c(H_1),...,l^D_c(H_t)\}=l^D_c(H_t)$.
We construct a sequence of graphs by the following steps:

Step 0: \emph{$L=1$ and $i=0$.}

Step 1: \emph{$G_{i+1}=H_{L,L-1+l^D_c(H_L)}$. $L=L+l^D_c(H_L)$ and $i=i+1$.}

Step 2: If $H_t$ is a subgraph of $G_i$, then the sequence is completed.

Step 3: \emph{Go to Step 1.}

Since $l^D_c(H_1)\leq t$, we can obtain a sequence of graphs at certain integer $i$ after step 2.
If $l^D_c(H_1)=t$, then the construction is completed with $i=1$
and $G_1=H_{1,t}$. Hence $\{G_1\}$ is a decomposition of $G$ and $t$ is the smallest positive integer $h_i$ such that
$f_D(H_{1,h_i})\geq ch_i$.

Suppose that $l^D_c(H_1)<t$. Then the construction is completed with $i>1$.
Let $G_k$ be the last graph of the sequence.
By the maximality of $l^D_c(H_1)$, $H_1$ is not a subgraph of $G_k$.
Hence $\{G_1,...,G_k\}$ is a decomposition of $G$.
By the construction of $G_i$,
we have $G_i=H_{r_i,r_i-1+h_i}, 1\leq i\leq k$,
such that $h_i$ is the smallest positive integer satisfying
$f_D(H_{r_i,r_i-1+h_i})\geq ch_i$.
\end{proof}

\begin{thm}\label{cor}
Suppose that a graph $G$ is generalized periodic with a transitive decomposition $\{H_1,...,H_t\}$.
Let $c$ be a positive real number and let $\Sigma$ be an arbitrary surface.
Then $\textrm{cr}_{\Sigma}(G)\geq\lceil ct\rceil$ if and only if
max$\{l^D_c(H_1),...,l^D_c(H_t)\}\leq t$ for each good drawing $D$ of $G$ in $\Sigma$.
\end{thm}

\begin{proof} Assume that $1\leq i\leq t$ if not specified.
If $\textrm{cr}_{\Sigma}(G)\geq\lceil ct\rceil$, then $\textrm{cr}(D)\geq ct$ by the definition of crossing number.
By \eqref{I1} and \eqref{I2}, we have $f_D(H_{i,i-1+t})=\textrm{cr}(D)$.
It follows that $f_D(H_{i,i-1+t})\geq ct$, and then $l^D_i\leq t$.
Hence max$\{l^D_c(H_1),...,$ $l^D_c(H_t)\}\leq t$.

Suppose that max$\{l^D_c(H_1),...,l^D_c(H_t)\}\leq t$.
By Lemma \ref{decom}, there exists a decomposition $\{G_1,...,G_k\}$ of $G$ such that $G_i=H_{r_i,r_i-1+h_i}$ and $f_D(H_{r_i,r_i-1+h_i})\geq ch_i, 1\leq i\leq k$.
Hence $\sum_{i=1}^kf_D(G_i)\geq ct$.
By \eqref{I1} and \eqref{I2}, we have $\textrm{cr}(D)=\sum_{i=1}^kf_D(G_i)$.
Thus $\textrm{cr}(D)\geq ct$.
Since $\textrm{cr}(D)$ is an integer, we have $\textrm{cr}(D)\geq \lceil ct\rceil$.
It follows that $\textrm{cr}_{\Sigma}(G)\geq\lceil ct\rceil$ by the definition of crossing number.
\end{proof}

By Theorem \ref{cor}, we have the following result.
\begin{cor}\label{subcor}
Suppose that a graph $G$ is generalized periodic with a transitive decomposition $\{H_1,...,H_t\}$.
Let $c$ be a positive real number and let $\Sigma$ be an arbitrary surface.
If $D$ is a good drawing of $G$ in $\Sigma$ with $\textrm{cr}(D)<\lceil ct\rceil$,
then max$\{l^D_c(H_1),...,l^D_c(H_t)\}=t+1$.
\end{cor}

In the following sections, Corollary \ref{subcor} is applied to obtain the crossing number of a graph in the plane.
Suppose that $G$ is a graph, $H$ is a subgraph of $G$ and $D$ is a good drawing of $G$.
Let $D(H)$ denote the subdrawing of $H$ in $D$.
If an edge is not crossed by any other edges, we say that it is
\emph{clean} in $D$; otherwise, we say that it is \emph{crossed} in $D$.
We call $H$ is \emph{clean} in $D$ if each edge in $H$ is clean in $D$.
By the Jordan Curve Theorem, we have the following lemma.

\begin{lemma}\label{J}
In a graph $G$, let $C$ and $C'$ be two vertex-disjoint cycles and let $P_k=u_1u_2...u_k$ be a path of order $k$ with $V(P_k)\cap V(C)=\emptyset$.
Suppose that $D$ is a good drawing of $G$ in the plane.
Then $\textrm{cr}_D(C,C')$ is even; $\textrm{cr}_D(C,P_k)$ is even if $u_1$ and $u_k$ are in the same region of $D(C)$, and it is odd otherwise.
\end{lemma}

\section{The crossing number of the generalized Petersen graph $P(4k+2,2k)$ in the plane}\label{P4k2}
Let $V(P(4k+2,2k))=\{u_i,v_i:0\leq i\leq 4k+1\}$ and let
$E(P(4k+2,2k))=\{u_iu_{i+1},u_iv_i,v_iv_{i+2k}:0\leq i\leq 4k+1\}$ with the subscripts taken modulo $4k+2$, $k\geq 1$.
It is trivial that $P(6,2)$ is a planar graph.
Sara\v{z}in \cite{Sar} proved that $\textrm{cr}(P(10,4))=4$.
Lin et al. \cite{Lin09} showed the following results.
\begin{lemma}[\cite{Lin09}]\label{base}
$\textrm{cr}(P(14,6))=7$.
\end{lemma}
\begin{lemma}[\cite{Lin09}]\label{up}
$\textrm{cr}(P(4k+2,2k))\leq 2k+1, k\geq 3$.
\end{lemma}

Lin et al. \cite{Lin09} also conjectured that the equality holds in Lemma \ref{up}.
Our main result in this section is the following result, which confirms this conjecture.
\begin{thm}\label{thm1}
$\textrm{cr}(P(4k+2,2k))=2k+1$, $k\geq 3$.
\end{thm}

Let $E_i$ be a subgraph of $P(4k+2,2k)$ with $V(E_i)=\{u_i,u_{i+1},v_i,v_{i+2k+2}\}$ and
$E(E_i)=\{u_iu_{i+1},u_iv_i,v_iv_{i+2k+2}\}$, $0\leq i\leq 4k+1$.
Assume that $F_i=E_i\cup E_{i+2k+1}$, $0\leq i\leq 2k$.
It is easy to check that $\{F_0,...,F_{2k}\}$ is a transitive decomposition of $P(4k+2,2k)$.
Let $D$ be a good drawing of $P(4k+2,2k)$.
For $0\leq i\leq 2k$, let $f_D(F_i)$ be a function to count the number of crossings on $F_i$ in $D$ as follows.
\begin{equation}\label{p1}
f_D(F_i)=\textrm{cr}_D(F_i)+\textrm{cr}_D(F_i,P(4k+2,2k)\setminus E(F_i))/2.
\end{equation}
Let $V(B_i)=\{u_i,v_i\}$ and let $E(B_i)=\{u_iv_i\}$, $0\leq i\leq 4k+1$.
For brevity, let $B_{i,i+1}=B_i\cup B_{i+1}, 0\leq i\leq 4k+1$.
The following lemma is based on the structure of $P(4k+2,2k)$.

\begin{lemma}\label{condition}
Let $\textrm{cr}(P(4k-2,2(k-1)))\geq 2k-1$ and let $D$ be a good drawing of $P(4k+2,2k)$ with $\textrm{cr}(D)<2k+1$, $k\geq 4$.
Then we have the following results.

$(i). \textrm{cr}_D(B_{i,i+1}\cup B_{i+2k+1,i+2k+2},P(4k+2,2k)\setminus E(B_{i,i+1}\cup B_{i+2k+1,i+2k+2}))+\textrm{cr}_D(B_{i,i+1}\cup B_{i+2k+1,i+2k+2})\leq 1$;

$(ii). \textrm{cr}_D(B_{i,i+1}\cup B_{i+2k+2,i+2k+3},P(4k+2,2k)\setminus E(B_{i,i+1}\cup B_{i+2k+2,i+2k+3}))+\textrm{cr}_D(B_{i,i+1}\cup B_{i+2k+2,i+2k+3})\leq 1, 0\leq i\leq 2k$.
\end{lemma}
\begin{proof}
It is easy to check that $P(4k+2,2k)\setminus E(B_{i,i+1}\cup B_{i+2k+1,i+2k+2})$ is a subdivision of $P(4k-2,2(k-1))$.
Since $\textrm{cr}(P(4k-2,2(k-1)))\geq 2k-1$ and $\textrm{cr}(D)<2k+1$,
a good drawing of a subdivision of $P(4k-2,2(k-1))$ with at least $2k-1$ crossings can be obtained by deleting $E(B_{i,i+1}\cup B_{i+2k+1,i+2k+2})$ from $D$.
Hence $\textrm{cr}_D(B_{i,i+1}\cup B_{i+2k+1,i+2k+2},P(4k+2,2k)\setminus E(B_{i,i+1}\cup B_{i+2k+1,i+2k+2}))+\textrm{cr}_D(B_{i,i+1}\cup B_{i+2k+1,i+2k+2})\leq 1$.
Similarly, we have $\textrm{cr}_D(B_{i,i+1}\cup B_{i+2k+2,i+2k+3},P(4k+2,2k)\setminus E(B_{i,i+1}\cup B_{i+2k+2,i+2k+3}))+\textrm{cr}_D(B_{i,i+1}\cup B_{i+2k+2,i+2k+3})\leq 1$.
\end{proof}
Given that $\textrm{cr}(P(4k-2,2(k-1)))\geq 2k-1$ and $D$ is a good drawing of $P(4k+2,2k)$ with $\textrm{cr}(D)<2k+1$, $k\geq 4$,
the following two lemmas describe properties of $D$.
In this section, let $H_i$ be the hexagon $u_iu_{i+1}u_{i+2}v_{i+2}v_{i+2k+2}v_iu_i$, $0\leq i\leq 4k+1$,
and let $O_j$ be the octagon $F_j\cup B_{j+1}\cup B_{j+2k+2}$, $0\leq j\leq 2k$.
\begin{lemma}\label{twoc}
Let $\textrm{cr}(P(4k-2,2(k-1)))\geq 2k-1$ and let $D$ be a good drawing of $P(4k+2,2k)$ with $\textrm{cr}(D)<2k+1$, $k\geq 4$.
If $l^D_1(F_i)>2$, then $\textrm{cr}_D(H_i,H_{i+2k+1})=0$, $0\leq i\leq 2k$.
\end{lemma}
\begin{proof} To the contrary, without loss of generality, suppose that $l^D_1(F_0)>2$ and $\textrm{cr}_D(H_0,H_{2k+1})\geq 1$.
By $l^D_1(F_0)>2$, we have $f_D(F_0)<1$ and $f_D(F_{0,1})<2$.
Hence $\textrm{cr}_D(F_{0,1})\leq 1$.
Since $\textrm{cr}_D(H_0,H_{2k+1})\geq 1$, $\textrm{cr}_D(H_0,H_{2k+1})\geq 2$ by Lemma \ref{J}.
Then $B_2\cup B_{2k+3}$ is crossed in $D(H_0\cup H_{2k+1})$ by $f_D(F_{0,1})<2$.
Without loss of generality, assume that $B_2$ is crossed in $D(H_0\cup H_{2k+1})$.
Then $B_2\cup B_{2k+3}$ cannot be crossed by any other edge in $D$ by Lemma \ref{condition}.
Moreover, by Lemma \ref{condition} and $f_D(F_{0,1})<2$, $\textrm{cr}_D(H_0,H_{2k+1})=2$ and $\textrm{cr}_D(F_{0,1})=1$.
We distinguish two cases:

Case 1: $\textrm{cr}_D(B_2,B_{2k+3})=1$.

Hence, by Lemma \ref{J}, $u_2$ and $v_2$ are in different regions in $D(H_{2k+1})$, and then $H_2\setminus B_2$ crosses $H_{2k+1}$.
Similarly, $H_{2k+3}\setminus B_{2k+3}$ crosses $H_0$.
Then $f_D(F_{0,1})\geq 2$, a contradiction.

Case 2: $\textrm{cr}_D(B_2,B_{2k+3})=0$.

Then $\textrm{cr}_D(B_2,H_{2k+1}\setminus B_{2k+3})=1$.
Hence, by Lemma \ref{J}, $u_2$ and $v_2$ are in different regions in $D(H_{2k+1})$, and then $H_2\setminus B_2$ crosses $H_{2k+1}$.
It follows that $f_D(F_{0,1})\geq 2$, a contradiction.
\end{proof}

\begin{lemma}\label{2E}
Let $\textrm{cr}(P(4k-2,2(k-1)))\geq 2k-1$ and let $D$ be a good drawing of $P(4k+2,2k)$ with $\textrm{cr}(D)<2k+1$, $k\geq 4$.
If $l^D_1(F_i)>2$ and $\textrm{cr}_D(F_{i,i+1})=1$, then each of $u_{i+2}$, $v_{i+2}$,
$u_{i+2k+3}$ and $v_{i+2k+3}$ is in the same region in $D(O_i)$; and both $B_{i+2}$ and $B_{i+2k+3}$ are clean in $D(F_{i,i+1}\cup B_{i+2}\cup B_{i+2k+3})$, $0\leq i\leq 2k$.
\end{lemma}
\begin{proof} Without loss of generality, suppose that $l^D_1(F_0)>2$ and $\textrm{cr}_D(F_{0,1})=1$.
By $l^D_1(F_0)>2$, we have $f_D(F_0)<1$ and $f_D(F_{0,1})<2$.
If $u_2$ and $v_2$ are in different regions in $D(O_0)$, then
both $B_2$ and $H_2\setminus B_2$ cross $O_0$ in $D$ by Lemma \ref{J}.
Hence $f_D(F_{0,1})\geq 2$, a contradiction.
Then $u_2$ and $v_2$ are in the same region in $D(O_0)$.
Similarly, $u_{2k+3}$ and $v_{2k+3}$ are in the same region in $D(O_0)$.
If $u_2$ and $v_{2k+3}$ are in different regions in $D(O_0)$,
then both $u_2u_3v_3v_{2k+3}$ and $u_{2k+3}u_{2k+4}v_{2k+4}v_2$ cross $O_0$ by Lemma \ref{J}.
We have $f_D(F_{0,1})\geq 2$, a contradiction.
Thus each of $u_2, v_2, u_{2k+3}$ and $v_{2k+3}$ is in the same region in $D(O_0)$.

By $f_D(F_{0,1})<2$ and Lemma \ref{twoc}, we have $\textrm{cr}_D(H_0,H_{2k+1})=0$.
Since $D$ is a good drawing,
if $B_2$ is crossed in $D(F_{0,1}\cup B_2\cup B_{2k+3})$,
then $B_2$ crosses $O_0$ by $\textrm{cr}_D(H_0,H_{2k+1})=0$.
Since $u_2$ and $v_2$ are in the same region in $D(O_0)$,
$\textrm{cr}_D(B_2,O_0)\geq 2$ by Lemma \ref{J}, which contradicts Lemma \ref{condition}.
Then $B_2$ is clean in $D(F_{0,1}\cup B_2\cup B_{2k+3})$.
Similarly, $B_{2k+3}$ is clean in $D(F_{0,1}\cup B_2\cup B_{2k+3})$.
\end{proof}
Given that $\textrm{cr}(P(4k-2,2(k-1)))\geq 2k-1$ and $D$ is a good drawing of $P(4k+2,2k)$ with $\textrm{cr}(D)<2k+1$, $k\geq 4$, the following four lemmas show that
if $f_D(F_i)<1$ and $\textrm{cr}_D(F_{i,i+1})\geq 1$, then $f_D(F_{i,i+1})\geq 2$ or $f_D(F_{i,i+2})\geq 3$, $0\leq i\leq 2k$.
For a path $P$, if $u,v\in V(P)$, then $uPv$ denotes the consecutive vertices of $P$ from $u$ to $v$.
In this section, let $P_0$, $P_1$ and $P_2$ be the paths $u_1...u_{4k+1}u_0$, $v_1v_{2k+3}v_3v_{2k+5}...v_{4k+1}v_{2k+1}$ and $v_{2k+2}v_2v_{2k+4}v_4...v_{2k}v_0$ in $P(4k+2,2k)$, respectively.
\begin{lemma}\label{ocross}
Let $\textrm{cr}(P(4k-2,2(k-1)))\geq 2k-1$ and $D$ let be a good drawing of $P(4k+2,2k)$ with $\textrm{cr}(D)<2k+1$, $k\geq 4$.
If $l^D_1(F_i)>2$, then $\textrm{cr}_D(O_i)=0$, $0\leq i\leq 2k$.
\end{lemma}
\begin{proof}To the contrary, without loss of generality, suppose that $l^D_1(F_0)>2$ and $\textrm{cr}_D(O_0)\geq 1$.
By $l^D_1(F_0)>2$, we have $f_D(F_0)<1$ and $f_D(F_{0,1})<2$.
By $f_D(F_{0,1})<2$, we have $\textrm{cr}_D(O_0)=1$.
By $f_D(F_0)<1$ and $\textrm{cr}_D(O_0)=1$, $B_1\cup B_{2k+2}$ is crossed in $D(O_0)$.
Without loss of generality, assume that $B_1$ is crossed in $D(O_0)$.
Hence $B_2\cup B_{2k+3}$ is clean in $D$ by Lemma \ref{condition}, and $D(F_{0,1}\cup B_2\cup B_{2k+3})$ is isomorphic to one of the drawings in Fig. \ref{C1} by Lemma \ref{2E}.

If $D(F_{0,1}\cup B_2\cup B_{2k+3})$ is isomorphic to one of the drawings in Fig. \ref{C1}a, \ref{C1}b and \ref{C1}c,
then both $v_2P_2v_{2k}u_{2k}u_{2k+1}$ and $u_2u_3v_3P_1v_{2k+1}$ cross $F_{0,1}\cup B_2\cup B_{2k+3}$ in $D$ by Lemma \ref{J}.
Otherwise, $D(F_{0,1}\cup B_2\cup B_{2k+3})$ is isomorphic to the drawing in Fig. \ref{C1}d or \ref{C1}e.
Then both $u_{2k+3}u_{2k+4}v_{2k+4}P_2v_0$ and $v_{2k+3}P_1v_{4k+1}u_{4k+1}u_0$ cross $F_{0,1}\cup B_2\cup B_{2k+3}$ in $D$ by Lemma \ref{J}.
In $D$, since $B_2\cup B_{2k+3}$ is clean,
both paths cross $F_{0,1}$ in each cases.
Hence $f_D(F_{0,1})\geq 2$, a contradiction.
\end{proof}
\begin{lemma}\label{uvcross}
Let $\textrm{cr}(P(4k-2,2(k-1)))\geq 2k-1$ and let $D$ be a good drawing of $P(4k+2,2k)$ with $\textrm{cr}(D)<2k+1$, $k\geq 4$.
If $l^D_1(F_i)>2$ and $\textrm{cr}_D(O_i)=0$, then $\textrm{cr}_D(u_{i+1}u_{i+2}\cup u_{i+2k+2}u_{i+2k+3}\cup v_{i+1}v_{i+2k+3}\cup v_{i+2k+2}v_{i+2})=0$,
$0\leq i\leq 2k$.
\end{lemma}
\begin{proof}
To the contrary, without loss of generality,
suppose that $\textrm{cr}_D(u_1u_2\cup u_{2k+2}u_{2k+3}\cup v_1v_{2k+3}\cup v_{2k+2}v_2)\geq 1$ if $l^D_1(F_0)>2$ and $\textrm{cr}_D(O_0)=0$.
By $l^D_1(F_0)>2$, we have $f_D(F_0)<1$ and $f_D(F_{0,1})<2$.
By Lemma \ref{twoc}, we have $\textrm{cr}_D(u_1u_2\cup v_{2k+2}v_2, u_{2k+2}u_{2k+3}\cup v_1v_{2k+3})=0$.
Hence $\textrm{cr}_D(u_{2k+2}u_{2k+3}, v_1v_{2k+3})+\textrm{cr}_D(u_1u_2,$ $v_{2k+2}v_2)=1$ by $f_D(F_{0,1})<2$.
Without loss of generality, assume that $\textrm{cr}_D(u_1u_2,$ $v_{2k+2}v_2)=1$.
Then $D(F_{0,1}\cup B_2\cup B_{2k+3})$ is isomorphic to one of the drawings in Fig. \ref{C03} by Lemma \ref{2E}.
By Lemma \ref{J}, both $u_{2k+3}P_0u_0$ and $v_{2k+3}v_3u_3P_0u_{2k}v_{2k}v_0$ cross $v_2v_{2k+2}u_{2k+2}u_{2k+1}v_{2k+1}v_1u_1u_2$.
Hence $f_D(F_{0,1})\geq 2$, a contradiction.
\end{proof}
\begin{lemma}\label{vcross}
Let $\textrm{cr}(P(4k-2,2(k-1)))\geq 2k-1$ and let $D$ be a good drawing of $P(4k+2,2k)$ with $\textrm{cr}(D)<2k+1$, $k\geq 4$.
If $l^D_1(F_i)>3$, $\textrm{cr}_D(O_i)=0$ and $\textrm{cr}_D(v_{i+1}v_{i+2k+3}\cup v_{i+2k+2}v_{i+2},O_i)=0$,
$0\leq i\leq 2k$.
\end{lemma}
\begin{proof}
To the contrary, without loss of generality,
suppose that $\textrm{cr}_D(v_1v_{2k+3}\cup v_{2k+2}v_2,O_0)\geq 1$ if $l^D_1(F_0)>3$ and $\textrm{cr}_D(O_0)=0$.
By $l^D_1(F_0)>3$, we have $f_D(F_0)<1$, $f_D(F_{0,1})<2$ and $f_D(F_{0,2})<3$.
Since $\textrm{cr}_D(v_1v_{2k+3}\cup v_{2k+2}v_2,O_0)\geq 1$, we have $\textrm{cr}_D(v_1v_{2k+3}\cup v_{2k+2}v_2,O_0)=1$ by $f_D(F_{0,1})<2$.
Without loss of generality, assume that $\textrm{cr}_D(v_{2k+2}v_2,O_0)=1$.
By Lemma \ref{twoc}, $\textrm{cr}_D(v_{2k+2}v_2,E_{2k+1})=0$.
Then $D(F_{0,1}\cup B_2\cup B_{2k+3})$ is isomorphic to one of the drawings in Fig. \ref{E2C1}a, \ref{E2C1}b and \ref{E2C1}c by Lemma \ref{2E}, since $D$ is a good drawing.
\captionsetup[figure]{labelfont={bf},name={Fig.},labelsep=period}
\begin{figure}
  \centering
  \subfigure[]{
    \label{fig:subfig:a} 
    \includegraphics[width=0.17\textwidth]{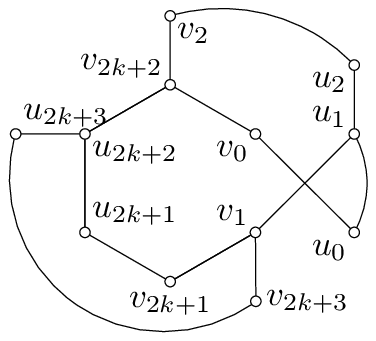}}
  \subfigure[]{
    \label{fig:subfig:a} 
    \includegraphics[width=0.17\textwidth]{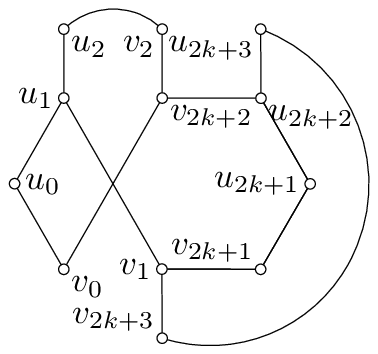}}
  \subfigure[]{
    \label{fig:subfig:b} 
    \includegraphics[width=0.22\textwidth]{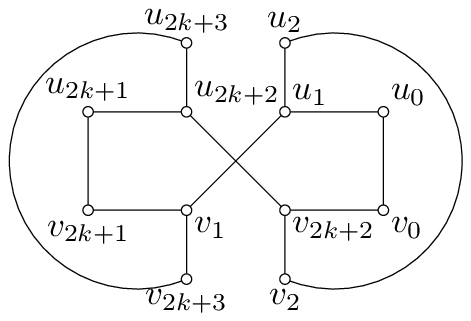}}
    \subfigure[]{
    \label{fig:subfig:b} 
    \includegraphics[width=0.17\textwidth]{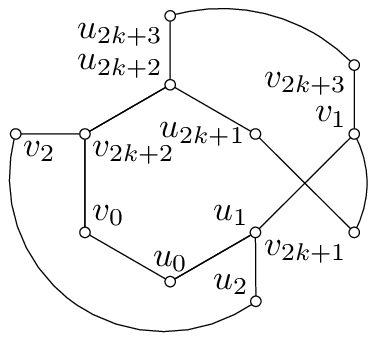}}
    \subfigure[]{
    \label{fig:subfig:b} 
    \includegraphics[width=0.17\textwidth]{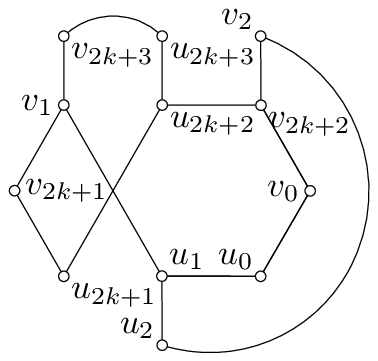}}
  \caption{Good drawings of $D(F_{0,1}\cup B_2\cup B_{2k+3})$ with $\textrm{cr}_D(O_0)=1$.}\label{C1}
\end{figure}
\begin{figure}
\centering
\subfigure[]{
\label{fig:subfig:a} 
\includegraphics[width=0.23\textwidth]{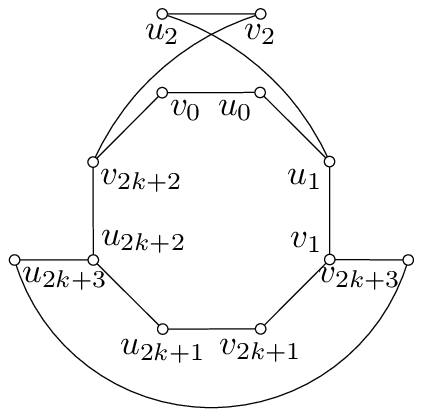}}
\subfigure[]{
\label{fig:subfig:a} 
\includegraphics[width=0.23\textwidth]{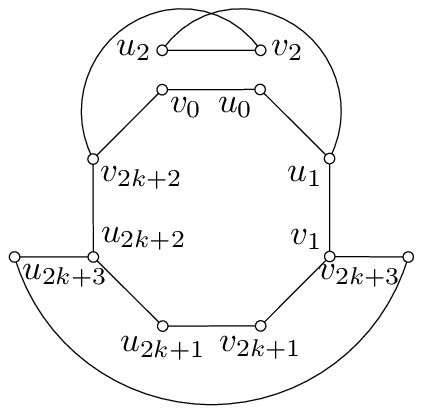}}
\caption{Good drawings of $D(F_{0,1}\cup B_2\cup B_{2k+3})$ with $\textrm{cr}_D(u_1u_2,v_{2k+2}v_2)=1$.}\label{C03}
\end{figure}

If $D(F_{0,1}\cup B_2\cup B_{2k+3})$ is isomorphic to the drawing in Fig. \ref{E2C1}a,
then both $v_0v_{2k}u_{2k}u_{2k+1}$ and $u_0u_{4k+1}v_{4k+1}v_{2k+1}$ cross $F_{0,1}\cup B_2\cup B_{2k+3}$ by Lemma \ref{J}.
Since $B_1$ is crossed, $B_2\cup B_{2k+3}$ is clean by Lemma \ref{condition}.
Thus $f_D(F_{0,1})\geq 2$, a contradiction.
If $D(F_{0,1}\cup B_2\cup B_{2k+3})$ is isomorphic to the drawing in Fig. \ref{E2C1}b,
then both $v_2P_2v_{2k}u_{2k}u_{2k+1}$ and $u_2u_3v_3P_1v_{2k+1}$ cross $u_{2k+3}v_{2k+3}v_1u_1u_0v_0v_{2k+2}u_{2k+2}u_{2k+3}$ by Lemma \ref{J}.
Since $B_0$ is crossed, $B_{2k+3}$ is clean in $D$ by Lemma \ref{condition}.
Hence $f_D(F_{0,1})\geq 2$, a contradiction.

Suppose that $D(F_{0,1}\cup B_2\cup B_{2k+3})$ is isomorphic to the drawing in Fig. \ref{E2C1}c.
By $f_D(F_{0,1})<2$ and $\textrm{cr}_D(F_{0,1})=1$,
$F_{0,1}$ is crossed at most once by any other edge.
By Lemma \ref{condition},
$B_2\cup B_{2k+3}$ is crossed at most once.
By Lemma \ref{J}, both $v_0v_{2k}u_{2k}u_{2k+1}$ and $u_0u_{4k+1}v_{4k+1}v_{2k+1}$ cross $v_2v_{2k+2}u_{2k+2}u_{2k+3}v_{2k+3}v_1u_1u_0$.
It follows that one of the two paths crosses $B_{2k+3}$ exactly once.
Since $B_{2k+3}$ is crossed, $B_3\cup B_{2k+4}$ is clean in $D$ by Lemma \ref{condition}.
Hence $D(F_{0,2}\cup B_3\cup B_{2k+4})$ is isomorphic to the drawing in Fig. \ref{E2C1}d by $f_D(F_{0,2})<3$.

\begin{figure}
  \centering
    \subfigure[]{
    \label{fig:subfig:a} 
    \includegraphics[width=0.23\textwidth]{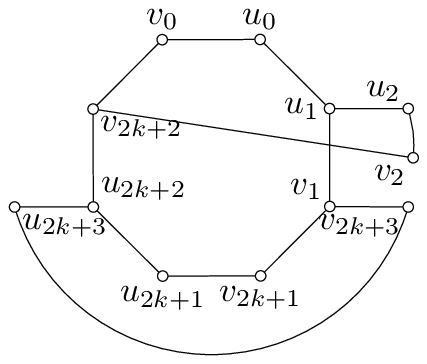}}
  \subfigure[]{
    \label{fig:subfig:a} 
    \includegraphics[width=0.23\textwidth]{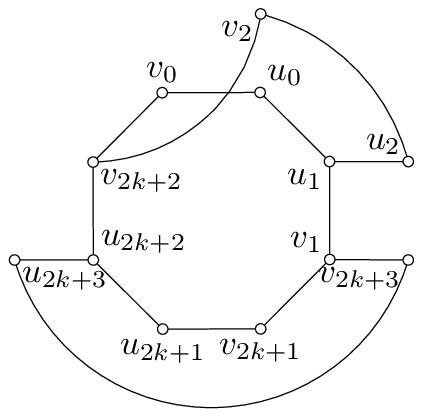}}
  \subfigure[]{
    \label{fig:subfig:b} 
    \includegraphics[width=0.23\textwidth]{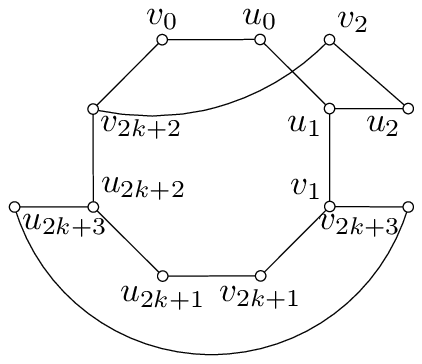}}
  \subfigure[]{
    \label{fig:subfig:b} 
    \includegraphics[width=0.23\textwidth]{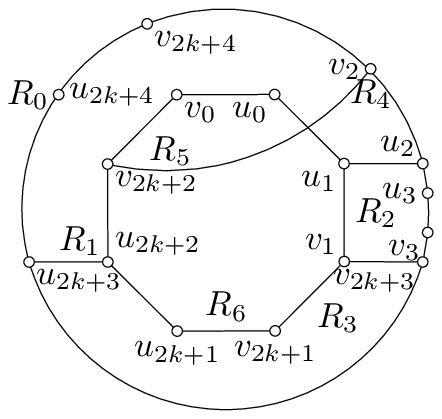}}
\caption{Subdrawings from $D$.}\label{E2C1}
\end{figure}

Without loss of generality, assume that $v_0v_{2k}u_{2k}u_{2k+1}$ crosses $B_{2k+3}$.
Then by Lemma \ref{J}, $v_0v_{2k}u_{2k}u_{2k+1}$ crosses $O_2$ at least twice.
Similarly, if $u_{4k+1}$ is in $R_0$, as shown in Fig. \ref{E2C1}d,
then $u_0u_{4k+1}v_{4k+1}v_{2k+1}$ crosses $O_2$ at least twice.
If $u_{4k+1}$ is in $R_2$,
then $u_0u_{4k+1}v_{4k+1}v_{2k+1}$ crosses $H_1$ at least twice by Lemma \ref{J}.
In both cases above, we have $f_D(F_{0,2})\geq 3$, a contradiction.
Then $u_{4k+1}$ is in one of the regions $R_1$, $R_3$, $R_4$, $R_5$ and $R_6$.
By Lemma \ref{J}, $v_3v_{2k+5}u_{2k+5}P_0u_{4k+1}$ crosses $F_{0,2}\cup B_3\cup B_{2k+4}$.
Then $f_D(F_{0,2})\geq 3$, a contradiction.
\end{proof}
Though $v_{2k+2}v_2$ and $u_{2k+2}u_{2k+3}$ are not in a common equivalence class on $E(P(4k+2,2k))$ with respect to automorphism of $P(4k+2,2k)$,
the proof of the following lemma is quite similar with the proof above.
\begin{lemma}\label{ucross}
Let $\textrm{cr}(P(4k-2,2(k-1)))\geq 2k-1$ and let $D$ be a good drawing of $P(4k+2,2k)$ with $\textrm{cr}(D)<2k+1$, $k\geq 4$.
If $l^D_1(F_i)>3$ and $\textrm{cr}_D(O_i)=0$, then $\textrm{cr}_D(u_{i+1}u_{i+2}\cup u_{i+2k+2}u_{i+2k+3},O_i)=0$, $0\leq i\leq 2k$.
\end{lemma}
\begin{proof}
To the contrary, without loss of generality,
suppose that $\textrm{cr}_D(u_1u_2\cup u_{2k+2}u_{2k+3},O_0)\geq 1$ if $l^D_1(F_0)>3$ and $\textrm{cr}_D(O_0)=0$.
By $l^D_1(F_0)>3$, we have $f_D(F_0)<1$, $f_D(F_{0,1})<2$ and $f_D(F_{0,2})<3$.
Since $\textrm{cr}_D(u_1u_2\cup u_{2k+2}u_{2k+3},O_0)\geq 1$, we have $\textrm{cr}_D(u_1u_2\cup u_{2k+2}u_{2k+3},O_0)=1$ by $f_D(F_{0,1})<2$.
Without loss of generality, assume that $\textrm{cr}_D(u_{2k+2}u_{2k+3},O_0)=1$.
By Lemma \ref{twoc}, $\textrm{cr}_D(u_{2k+2}u_{2k+3},$ $E_0)=0$.
Then $D(F_{0,1}\cup B_2\cup B_{2k+3})$ is isomorphic to one of the drawings in Fig. \ref{E2C2}a, \ref{E2C2}b and \ref{E2C2}c by Lemma \ref{2E}, since $D$ is a good drawing.

If $D(F_{0,1}\cup B_2\cup B_{2k+3})$ is isomorphic to the drawing in Fig. \ref{E2C2}a,
then $B_2\cup B_{2k+3}$ is clean by Lemma \ref{condition}.
By Lemma \ref{J}, both $v_0v_{2k}u_{2k}u_{2k+1}$ and $u_0u_{4k+1}v_{4k+1}v_{2k+1}$ cross $F_{0,1}\cup B_2\cup B_{2k+3}$.
Thus $f_D(F_{0,1})\geq 2$, a contradiction.
If $D(F_{0,1}\cup B_2\cup B_{2k+3})$ is isomorphic to the drawing in Fig. \ref{E2C2}b,
then $B_2$ is clean by Lemma \ref{condition}.
Both $u_{2k+3}P_0u_0$ and $v_{2k+3}v_3u_3P_0u_{2k}v_{2k}v_0$ cross $u_{2k+2}v_{2k+2}v_2u_2u_1v_1v_{2k+1}u_{2k+1}u_{2k+2}$ by Lemma \ref{J}.
Hence $f_D(F_{0,1})\geq 2$, a contradiction.

Suppose that $D(F_{0,1}\cup B_2\cup B_{2k+3})$ is isomorphic to the drawing in Fig. \ref{E2C2}c.
Since $\textrm{cr}_D(F_{0,1})=1$,
$F_{0,1}$ is crossed at most once by any other edge by $f_D(F_{0,1})<2$.
By Lemma \ref{condition},
$B_2$ is crossed at most once.
By Lemma \ref{J}, both $v_0v_{2k}u_{2k}u_{2k+1}$ and $u_0u_{4k+1}v_{4k+1}v_{2k+1}$ cross $u_{2k+3}u_{2k+2}v_{2k+2}v_2u_2u_1v_1v_{2k+1}$.
It follows that one of the two paths crosses $B_2$ exactly once.
Since $B_2$ is crossed, $B_3\cup B_{2k+4}$ is clean in $D$ by Lemma \ref{condition}.
Hence $D(F_{0,2}\cup B_3\cup B_{2k+4})$ is isomorphic to the drawing in Fig. \ref{E2C2}d by $f_D(F_{0,2})<3$.

Without loss of generality, assume that $v_0v_{2k}u_{2k}u_{2k+1}$ crosses $B_2$.
By Lemma \ref{J}, $v_0v_{2k}u_{2k}u_{2k+1}$ crosses $O_2$ at least twice.
Similarly, if $u_{4k+1}$ is in $R_0$, as shown in Fig. \ref{E2C2}d,
then $u_0u_{4k+1}v_{4k+1}v_{2k+1}$ crosses $O_2$ at least twice.
If $u_{4k+1}$ is in $R_2$,
then $u_0u_{4k+1}v_{4k+1}v_{2k+1}$ crosses $H_1$ at least twice by Lemma \ref{J}.
In both cases above, we have $f_D(F_{0,2})\geq 3$, a contradiction.
Then $u_{4k+1}$ is in one of the regions $R_1$, $R_3$, $R_4$, $R_5$ and $R_6$.
By Lemma \ref{J}, $v_3v_{2k+5}u_{2k+5}P_0u_{4k+1}$ crosses $F_{0,2}\cup B_3\cup B_{2k+4}$.
Then $f_D(F_{0,2})\geq 3$, a contradiction.
\end{proof}

  \begin{figure}
\centering
    \subfigure[]{
    \label{fig:subfig:a} 
    \includegraphics[width=0.23\textwidth]{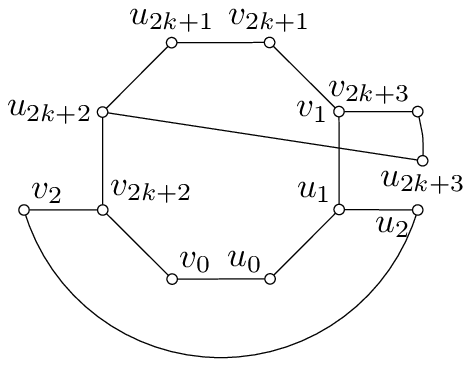}}
  \subfigure[]{
    \label{fig:subfig:a} 
    \includegraphics[width=0.23\textwidth]{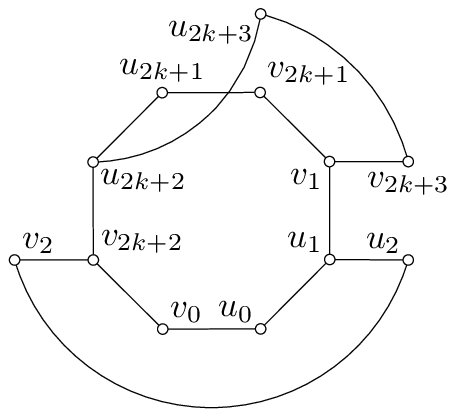}}
  \subfigure[]{
    \label{fig:subfig:b} 
    \includegraphics[width=0.23\textwidth]{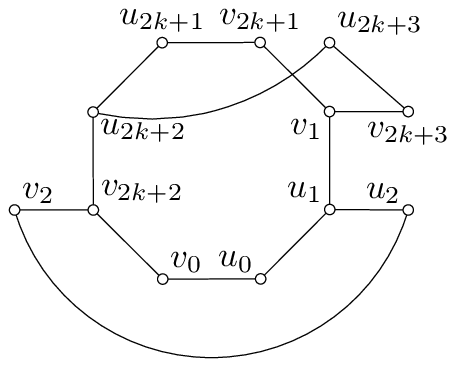}}
    \subfigure[]{
    \label{fig:subfig:b} 
    \includegraphics[width=0.23\textwidth]{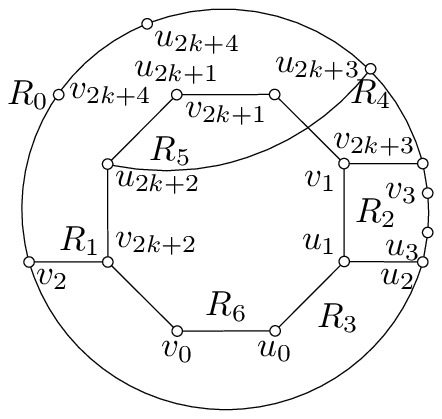}}
  \caption{Subdrawings from $D$.}\label{E2C2}
\end{figure}
Given that $\textrm{cr}(P(4k-2,2(k-1)))\geq 2k-1$ and $D$ is a good drawing of $P(4k+2,2k)$ with $\textrm{cr}(D)<2k+1$, $k\geq 4$,
the following three lemmas show that if $f_D(F_i)<1$, $f_D(F_{i,i+1})<2$ and $\textrm{cr}_D(F_{i,i+1})=0$, then $f_D(F_{i,i+2})\geq 3$ or $f_D(F_{i,i+3})\geq 4$, $0\leq i\leq 2k$.
\begin{lemma}\label{drawingd}
Let $\textrm{cr}(P(4k-2,2(k-1)))\geq 2k-1$ and let $D$ be a good drawing of $P(4k+2,2k)$ with $\textrm{cr}(D)<2k+1$, $k\geq 4$.
For $0\leq i\leq 2k$, if $l^D_1(F_i)>2$ and $\textrm{cr}_D(F_{i,i+1})= 0$,
then $D(F_{i,i+1}\cup B_{i+2}\cup B_{i+2k+3})$ is isomorphic to the drawing in Fig. \ref{C02}d.
\end{lemma}
\begin{proof}Without loss of generality, assume that $l^D_1(F_0)>2$ and $\textrm{cr}_D(F_{0,1})=0$.
By $l^D_1(F_0)>2$, we have $f_D(F_0)<1$ and $f_D(F_{0,1})<2$.
We claim that $u_{2k+3}$ and $v_{2k+3}$ are in the same region in $D(O_0)$.
Suppose not.
Then both $B_{2k+3}$ and $H_{2k+3}\setminus B_{2k+3}$ cross $O_0$ by Lemma \ref{J}.
Hence $B_{2k+3}$ is crossed exactly once, and $B_2$ is clean by Lemma \ref{condition}.
If $\textrm{cr}_D(B_{2k+3},B_1\cup B_{2k+2})=0$,
then $B_1\cup B_{2k+2}$ is clean by Lemma \ref{condition}.
Hence $f_D(F_0)\geq 1$, a contradiction.
Otherwise, $\textrm{cr}_D(B_{2k+3},B_1\cup B_{2k+2})=1$.
Without loss of generality, assume that $\textrm{cr}_D(B_{2k+3},B_1)=1$.
Then $D(F_{0,1}\cup B_2\cup B_{2k+3})$ is isomorphic to the drawing in Fig. \ref{C02}a or \ref{C02}b.
By Lemma \ref{condition}, $B_2\cup B_{2k+3}$ cannot be crossed by any other edge.

By Lemma \ref{J}, if $D(F_{0,1}\cup B_2\cup B_{2k+3})$ is isomorphic to the drawing in Fig. \ref{C02}a,
then both $u_2P_0u_{2k+1}$ and $v_2P_2v_{4k}u_{4k}u_{4k+1}v_{4k+1}v_{2k+1}$ cross $F_{0,1}\cup B_2\cup B_{2k+3}$;
if $D(F_{0,1}\cup B_2\cup B_{2k+3})$ is isomorphic to the drawing in Fig. \ref{C02}b,
then both $u_0u_{4k+1}v_{4k+1}v_{2k+1}$ and $v_0v_{2k}u_{2k}u_{2k+1}$ cross $F_{0,1}\cup B_2\cup B_{2k+3}$.
Recall that $H_{2k+3}\setminus B_{2k+3}$ crosses $O_0$.
In both cases, $f_D(F_{0,1})\geq 2$, a contradiction.
Therefore, $u_{2k+3}$ and $v_{2k+3}$ are in the same region in $D(O_0)$.
Similarly, $u_2$ and $v_2$ are in the same region in $D(O_0)$.

Since $D$ is a good drawing, if $B_2$ is crossed in $D(F_{0,1}\cup B_2\cup B_{2k+3})$,
then $B_2$ crosses $O_0$ by Lemma \ref{twoc}.
By Lemma \ref{J}, $B_2$ crosses $O_0$ at least twice,
which contradicts Lemma \ref{condition}.
Hence $B_2$ is clean in $D(F_{0,1}\cup B_2\cup B_{2k+3})$.
Similarly, $B_{2k+3}$ is clean in $D(F_{0,1}\cup B_2\cup B_{2k+3})$.
Hence $D(F_{0,1}\cup B_2\cup B_{2k+3})$ is isomorphic to the drawing in Fig. \ref{C02}c or \ref{C02}d.

\begin{figure}
\centering
\subfigure[]{
    \label{fig:subfig:b} 
    \includegraphics[width=0.23\textwidth]{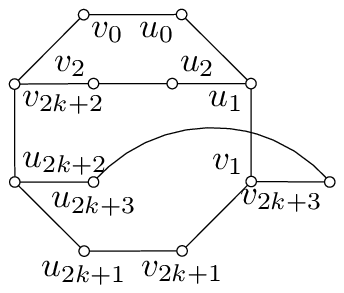}}
\subfigure[]{
    \label{fig:subfig:b} 
    \includegraphics[width=0.23\textwidth]{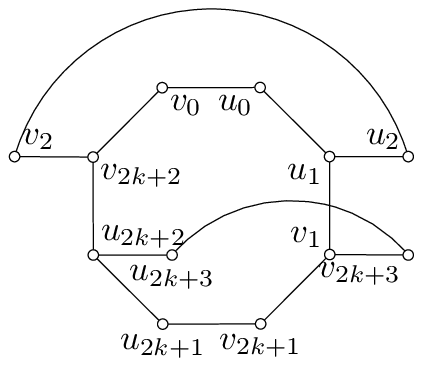}}
\subfigure[]{
    \label{fig:subfig:b} 
    \includegraphics[width=0.23\textwidth]{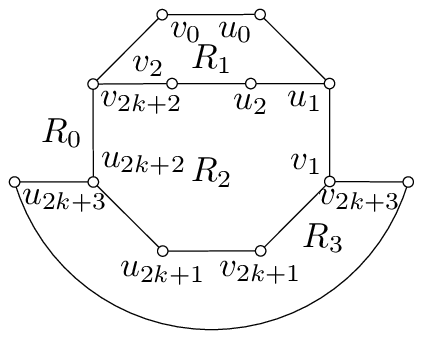}}
\subfigure[]{
    \label{fig:subfig:b} 
    \includegraphics[width=0.23\textwidth]{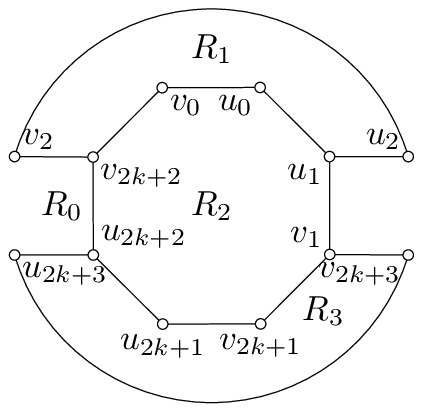}}
\caption{Good drawings of $F_{0,1}\cup B_2\cup B_{2k+3}$.}\label{C02}
\end{figure}
If $D(F_{0,1}\cup B_2\cup B_{2k+3})$ is isomorphic to the drawing in Fig. \ref{C02}c,
then both $u_2u_3v_3v_{2k+3}$ and $v_2v_{2k+4}u_{2k+4}u_{2k+3}$ cross $O_0$ by Lemma \ref{J}.
By Lemma \ref{condition},
$B_1\cup B_{2k+2}$ is crossed at most once.
It follows that one of the two paths crosses $B_1\cup B_{2k+2}$ by $f_D(F_0)<1$.
Then $B_2\cup B_{2k+3}$ is clean in $D$ by Lemma \ref{condition}.
By Lemma \ref{J}, both $v_0v_{2k}u_{2k}u_{2k+1}$ and $u_0u_{4k+1}v_{4k+1}v_{2k+1}$ cross $O_1$.
Hence $f_D(F_{0,1})\geq 2$, a contradiction.
Therefore, $D(F_{0,1}\cup B_2\cup B_{2k+3})$ is isomorphic to the drawing in Fig. \ref{C02}d.
\end{proof}

\begin{lemma}\label{drawinga}
Let $\textrm{cr}(P(4k-2,2(k-1)))\geq 2k-1$ and let $D$ be a good drawing of $P(4k+2,2k)$ with $\textrm{cr}(D)<2k+1$, $k\geq 4$.
For $0\leq i\leq 2k$, if $l^D_1(F_i)>3$ and $D(F_{i,i+1}\cup B_{i+2}\cup B_{i+2k+3})$ is isomorphic to the drawing in Fig. \ref{C02}d,
then $D(F_{i,i+2}\cup B_{i+3}\cup B_{i+2k+4})$ is isomorphic to the drawing in Fig. \ref{C02l}.
\end{lemma}
\begin{proof}Without loss of generality, assume that $l^D_1(F_0)>3$ and $D(F_{0,1}\cup B_2\cup B_{2k+3})$ is isomorphic to the drawing in Fig. \ref{C02}d.
Then $f_D(F_0)<1$, $f_D(F_{0,1})<2$ and $f_D(F_{0,2})<3$ by $l^D_1(F_0)>3$.
Let $C_1$ be $u_1v_1v_{2k+1}u_{2k+1}u_{2k+2}v_{2k+2}v_2u_2u_1$ and let $C_2$ be $u_1v_1v_{2k+3}u_{2k+3}u_{2k+2}v_{2k+2}v_0u_0u_1$.
Then $C_1\cap C_2=B_1\cup B_{2k+2}$.
Let $Q_1, Q_2, Q_3$ and $Q_4$ be
$v_{2k+3}P_1v_{2k-1}u_{2k-1}u_{2k}v_{2k}v_0$,
$u_{2k+3}P_0u_0$,
$u_{2}P_0u_{2k+1}$ and
$v_2P_2v_{4k}u_{4k}u_{4k+1}v_{4k+1}v_{2k+1}$, respectively.
It is easy to check that $Q_1\cap Q_3=u_{2k-1}u_{2k}$, $Q_2\cap Q_4=u_{4k}u_{4k+1}$ and $(Q_1\cup Q_3)\cap (Q_2\cup Q_4)=\emptyset$.
By Lemma \ref{J}, both $Q_1$ and $Q_2$ cross $C_1$, and both $Q_3$ and $Q_4$ cross $C_2$.

We claim that $\textrm{cr}_D(Q_1\cap Q_3,C_1\cap C_2)=0$.
Suppose not. Then $\textrm{cr}_D(u_{2k-1}u_{2k},$ $B_1\cup B_{2k+2})\geq 1$.
Then $\textrm{cr}_D(u_{2k-1}u_{2k},B_1\cup B_{2k+2})=1$, $B_1\cup B_{2k+2}$ cannot be crossed by any other edge and
$B_2\cup B_{2k+3}$ is clean by Lemma \ref{condition}.
Recall that $Q_2$ crosses $C_1$, that $Q_4$ crosses $C_2$, that $C_1\cap C_2=B_1\cup B_{2k+2}$ and
that $u_{2k-1}u_{2k}\cap (Q_2\cup Q_4)=\emptyset$.
Hence $\textrm{cr}_D(u_{2k-1}u_{2k}\cup Q_2\cup Q_4,C_1\cup C_2)\geq 3$.
Then $\textrm{cr}_D(Q_4,C_2)=1$ by $f_D(F_{0,1})<2$.
By $f_D(F_{0,1})<2$ and $E(v_{2k}v_{4k})\cap (E(Q_2\cup Q_4))=\emptyset$, we have $\textrm{cr}_D(v_{2k}v_{4k},O_0)=0$.
Similarly, $\textrm{cr}_D(v_{4k+1}v_{2k-1}\cup B_{2k}\cup u_2Q_3u_{2k-1},O_0)=0$.
Since $\textrm{cr}_D(u_2Q_3u_{2k-1},O_0)=0$, $u_{2k-1}$ is in the region $R_0$ or $R_1$, as shown in Fig. \ref{C02}d.
By $\textrm{cr}_D(u_{2k-1}u_{2k},B_1\cup B_{2k+2})=1$, $u_{2k}$ is in $R_2$ by Lemma \ref{J}.

Since $u_{2k}$ is in $R_2$ and $\textrm{cr}_D(B_{2k}\cup v_{2k}v_{4k},O_0)=0$, $v_{4k}$ is in $R_2$.
By Lemma \ref{J}, $v_2Q_4v_{4k}$ crosses $O_0$.
Hence, by $\textrm{cr}_D(Q_4,C_2)=1$, $\textrm{cr}_D(v_{4k}Q_4v_{4k+1},O_0)=0$ and then $v_{4k+1}$ is in $R_2$.
Moreover, $v_{2k-1}$ is in $R_2$ by $\textrm{cr}_D(v_{4k+1}v_{2k-1},O_0)=0$.
By Lemma \ref{J}, $\textrm{cr}_D(B_{2k-1},O_0)\geq 1$.
Thus $f_D(F_{0,1})\geq 2$ by $\textrm{cr}_D(u_{2k-1}u_{2k}\cup Q_2\cup Q_4,C_1\cup C_2)\geq 3$, a contradiction.
It follows that $\textrm{cr}_D(Q_1\cap Q_3,C_1\cap C_2)=0$.
Similarly, $\textrm{cr}_D(Q_2\cap Q_4,C_1\cap C_2)=0$.

Since $Q_1$ crosses $C_1$, $Q_3$ crosses $C_2$ and $\textrm{cr}_D(Q_1\cap Q_3,C_1\cap C_2)=0$,
we have $\textrm{cr}_D(Q_1\cup Q_3,C_1\cup C_2)\geq 2$.
Similarly, $\textrm{cr}_D(Q_2\cup Q_4,C_1\cup C_2)\geq 2$.
Since $(Q_1\cup Q_3)\cap (Q_2\cup Q_4)=\emptyset$, we have $\textrm{cr}_D(\cup_{i=1}^4 Q_i,C_1\cup C_2)\geq 4$.
Hence, by $f_D(F_{0,1})<2$, $\textrm{cr}_D(\cup_{i=1}^4 Q_i,B_2\cup B_{2k+3})\geq 1$.
Then $B_1\cup B_{2k+2}\cup B_3\cup B_{2k+4}$ is clean by Lemma \ref{condition}.
Moreover, by $f_D(F_{0,2})<3$ and $\textrm{cr}_D(\cup_{i=1}^4 Q_i,C_1\cup C_2)\geq 4$,
$\textrm{cr}_D(v_{2k+3}v_3\cup u_{2k+3}u_{2k+4}\cup u_2u_3\cup v_2v_{2k+4})=0$ and $\textrm{cr}_D(v_{2k+3}v_3\cup u_{2k+3}u_{2k+4}\cup u_2u_3\cup v_2v_{2k+4},C_1\cup C_2)\leq 1$.

We claim that $\textrm{cr}_D(v_{2k+3}v_3\cup u_{2k+3}u_{2k+4}\cup u_2u_3\cup v_2v_{2k+4},C_1\cup C_2)=0$.
Suppose not. Without loss of generality, suppose that $\textrm{cr}_D(u_2u_3,C_1\cup C_2)=1$.
Then by $f_D(F_{0,2})<3$,
$\textrm{cr}_D(v_{2k+3}v_3\cup u_{2k+3}u_{2k+4}\cup v_2v_{2k+4},C_1\cup C_2)=0$ and $\textrm{cr}_D(u_3Q_3u_{2k+1},C_1\cup C_2)=0$.
Hence $u_3$ is in $R_2$ or $R_3$.
Since $\textrm{cr}_D(v_{2k+3}v_3,C_1\cup C_2)=0$, $v_3$ is in $R_0$ or $R_3$.
Similarly, $u_{2k+4}$ is in $R_0$ or $R_3$, and $v_{2k+4}$ is in $R_0$ or $R_1$.
Since $B_3$ is clean, $v_3$ is in $R_3$.
Similarly, $u_{2k+4}$ is in $R_0$.
By $f_D(F_{0,2})<3$, $\textrm{cr}_D(v_3Q_1v_0,C_1)=1$.
Hence $v_{2k+5}$ is in $R_2$ or $R_3$ by Lemma \ref{J}.
By $f_D(F_{0,2})<3$, $\textrm{cr}_D(u_{2k+4}Q_2u_0,C_2)=1$.
Recall that $B_1\cup B_{2k+2}$ is clean in $D$.
Hence $u_{2k+5}$ is in $R_0$ or $R_1$ by Lemma \ref{J}.
Then $B_{2k+5}$ crosses $C_1\cup C_2$ by Lemma \ref{J}.
It follows that $f_D(F_{0,2})\geq 3$, a contradiction.
Thus $\textrm{cr}_D(v_{2k+3}v_3\cup u_{2k+3}u_{2k+4}\cup u_2u_3\cup v_2v_{2k+4},C_1\cup C_2)=0$.

Since $B_3\cup B_{2k+4}$ is clean, $\textrm{cr}_D(v_{2k+3}v_3\cup u_{2k+3}u_{2k+4}\cup u_2u_3\cup v_2v_{2k+4})=0$ and
$\textrm{cr}_D(v_{2k+3}v_3\cup u_{2k+3}u_{2k+4}\cup u_2u_3\cup v_2v_{2k+4},C_1\cup C_2)=0$,
$D(F_{0,2}\cup B_3\cup B_{2k+4})$ is isomorphic to the drawing in Fig. \ref{C02l}.
\end{proof}
\begin{lemma}\label{complete}
Let $\textrm{cr}(P(4k-2,2(k-1)))\geq 2k-1$ and let $D$ be a good drawing of $P(4k+2,2k)$ with $\textrm{cr}(D)<2k+1$, $k\geq 4$.
For $0\leq i\leq 2k$, if $l^D_1(F_i)>3$ and $D(F_{i,i+2}\cup B_{i+3}\cup B_{i+2k+4})$ is isomorphic to the drawing in Fig. \ref{C02l},
then $l^D_1(F_i)=4$.
\end{lemma}
\begin{proof}
To the contrary, without loss of generality, suppose that $l^D_1(F_0)>4$ if $D(F_{0,2}\cup B_3\cup B_{2k+4})$ is isomorphic to the drawing in Fig. \ref{C02l}.
By $l^D_1(F_0)>4$, we have $f_D(F_0)<1$, $f_D(F_{0,1})<2$, $f_D(F_{0,2})<3$ and $f_D(F_{0,3})<4$.
By Lemma \ref{J}, each of $u_3P_0u_{2k+1}$, $v_3P_1v_{2k+1}$, $u_{2k+4}P_0u_0$ and $v_{2k+4}P_2v_0$ crosses $O_1$.
By $f_D(F_{0,1})<2$, at least one of the four paths crosses $B_2\cup B_{2k+3}$.
Thus $B_1\cup B_{2k+2}\cup B_3\cup B_{2k+4}$ is clean by Lemma \ref{condition}.
By Lemma \ref{condition}, we have $\textrm{cr}_D(B_{2k}\cup B_{4k+1},F_{0,2}\cup B_3\cup B_{2k+4})\leq 1$.
Without loss of generality, suppose that $\textrm{cr}_D(B_{2k},F_{0,2}\cup B_3\cup B_{2k+4})=0$.
Then $u_{2k}$ and $v_{2k}$ are in the same region in $D(F_{0,2}\cup B_3\cup B_{2k+4})$.
We consider all the regions, as shown in Fig. \ref{C02l}, where $u_{2k}$ is in.
In each case, we obtain a contradiction.
\begin{figure}
\centering
\label{fig:subfig:b} 
\includegraphics[width=0.3\textwidth]{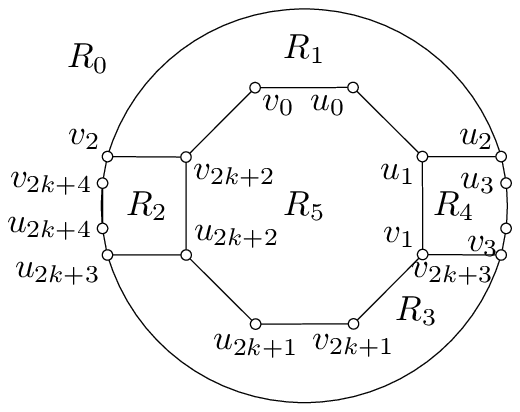}
\caption{A good drawing of $F_{0,2}\cup B_3\cup B_{2k+4}$.}\label{C02l}
\end{figure}

If $u_{2k}$ is in $R_5$,
then both $u_3P_0u_{2k}$ and $v_{2k+4}P_2v_{2k}$ cross $O_0$ by Lemma \ref{J}.
Since $B_1\cup B_{2k+2}$ is clean, we have $f_D(F_0)\geq 1$, a contradiction.
If $u_{2k}$ is in $R_1$,
then each of $u_{2k}u_{2k+1}$, $u_3P_0u_{2k}$ and $v_{2k+4}P_2v_{2k}$ crosses $H_0$ by Lemma \ref{J}.
Recall that both $u_{2k+4}P_0u_0$ and $v_3P_1v_{2k+1}$ cross $O_1$,
and that $B_2\cup B_{2k+3}$ is crossed at most once.
We have $f_D(F_{0,1})\geq 2$, a contradiction.
Similarly, $u_{2k}$ is not in $R_3$.

If $u_{2k}$ is in $R_0$,
then each of $v_{2k}v_{4k}u_{4k}u_{4k+1}u_0$, $u_{2k}u_{2k-1}v_{2k-1}v_{4k+1}v_{2k+1}$, $v_{2k}v_0$ and $u_{2k}u_{2k+1}$ crosses $O_2$ by Lemma \ref{J}.
If none of the four paths crosses $B_2\cup B_{2k+3}$, then each of them crosses $H_1\cup H_{2k+2}$ at least twice by Lemma \ref{J}.
Hence $f_D(F_{0,3})\geq 4$, a contradiction.
Otherwise, at least one of the four paths crosses $B_2\cup B_{2k+3}$.
By Lemma \ref{condition}, $B_3\cup B_{2k+4}$ is clean.
Recall that $B_2\cup B_{2k+3}$ is crossed at most once.
By the same argument, we have $f_D(F_{0,2})\geq 7/2$, a contradiction.

If $u_{2k}$ is in $R_2$, then each of $v_{2k}v_{4k}u_{4k}u_{4k+1}u_0$, $u_{2k}u_{2k-1}v_{2k-1}v_{4k+1}v_{2k+1}$,
$v_{2k}v_0$ and $u_{2k}u_{2k+1}$ crosses $O_1$ by Lemma \ref{J}.
If none of the four paths crosses $B_2\cup B_{2k+3}$, then $f_D(F_{0,1})\geq 2$, a contradiction.
Suppose that at least one of the four paths crosses $B_2\cup B_{2k+3}$.
By Lemma \ref{condition}, $B_1\cup B_{2k+1}\cup B_3\cup B_{2k+4}$ is clean,
and $B_2\cup B_{2k+3}$ is crossed exactly once.
By Lemma \ref{J}, the path crossing $B_2\cup B_{2k+3}$ crosses $O_2$ at least twice.
Hence $f_D(F_{0,2})\geq 5/2$.
By $f_D(F_{0,2})<3$, each of the three paths that not crossing $B_2\cup B_{2k+3}$ crosses $O_1$ only once.
Since $B_1\cup B_{2k+1}$ is clean, if $u_{2k}u_{2k-1}v_{2k-1}v_{4k+1}v_{2k+1}$ does not cross $B_2\cup B_{2k+3}$,
then $u_{2k-1}$ is in $R_2$ or $R_3$ by Lemma \ref{J}.
By Lemma \ref{J}, $u_3P_0u_{2k-1}$ crosses $F_{0,2}\cup B_3\cup B_{2k+4}$.
Since $B_3\cup B_{2k+4}$ is clean, we have $f_D(F_{0,2})\geq 3$, a contradiction.
Otherwise, $v_{2k}v_{4k}u_{4k}u_{4k+1}u_0$ does not cross $B_2\cup B_{2k+3}$.
By the same argument, $v_3P_1v_{4k+1}u_{4k+1}$ crosses $F_{0,2}\cup B_3\cup B_{2k+4}$ and
$f_D(F_{0,2})\geq 3$, a contradiction.
Similarly, $u_{2k}$ is not in $R_4$.
\end{proof}
We are ready to give the proof of Theorem \ref{thm1}.

\noindent\textit{Proof of Theorem \ref{thm1}.}
By Lemma \ref{up}, it suffices to prove $\textrm{cr}(P(4k+2,2k))\geq 2k+1$, $k\geq 3$.
We prove this inequality by induction on $k$.
By Lemma \ref{base}, $\textrm{cr}(P(4k+2,2k))\geq 2k+1$ if $k=3$.
Assume that $\textrm{cr}(P(4k-2,2(k-1)))\geq 2k-1$, $k\geq 4$.
We prove $\textrm{cr}(P(4k+2,2k))\geq 2k+1$ by contradiction.
Suppose that $\textrm{cr}(P(4k+2,2k))<2k+1$.
Then there exists a good drawing $D$ of $P(4k+2,2k)$ with $\textrm{cr}(D)<2k+1$.

Without loss of generality, assume that max$\{l^D_1(F_0),...,l^D_1(F_{2k})\}=l^D_1(F_0)$.
By Corollary \ref{subcor}, $l^D_1(F_0)=2k+2$.
By Lemma \ref{ocross}, Lemma \ref{uvcross}, Lemma \ref{vcross} and Lemma \ref{ucross},
we have $\textrm{cr}_D(F_{0,1})=0$.
Then $D(F_{0,1}\cup B_2\cup B_{2k+3})$ is isomorphic to the drawing in Fig. \ref{C02}d by Lemma \ref{drawingd}.
Hence $D(F_{0,2}\cup B_3\cup B_{2k+4})$ is isomorphic to the drawing in Fig. \ref{C02l} by Lemma \ref{drawinga}.
By Lemma \ref{complete}, we have $l^D_1(F_0)\leq 4$, a contradiction.$\hfill\Box$

\begin{remark}\label{rem1}
As shown in the proof above, Theorem \ref{cor} can be applied to construct algorithm to compute the lower bounds on the crossing number of generalized periodic graphs in an arbitrary surface.
\end{remark}
For example, we give algorithms to compute the lower bounds on the crossing number of generalized periodic graphs in the plane in the next section.

\section{Computing lower bounds on crossing numbers of generalized periodic graphs in the plane}

Suppose that $c$ is a positive number, and $\{H_1,...,H_{k}\}$ is a transitive decomposition of a graph $G$.
Let $\hat{\mathcal{D}}$ be the set of all good drawings of $G$ such that $\textrm{cr}(D)<ck$ for all $D\in\hat{\mathcal{D}}$.
By the definition of crossing number, $\textrm{cr}(G)<ck$ if and only if $\hat{\mathcal{D}}\not=\emptyset$,
and $\textrm{cr}(G)\geq ck$ if and only if $\hat{\mathcal{D}}=\emptyset$.
Suppose that $\textrm{cr}(G)<ck$ and
$\bar{D}$ is a good drawing of $G$ with $\textrm{cr}(\bar{D})<ck$.

By \eqref{I1}, \eqref{I2} and \eqref{fun}, we have $\textrm{cr}(\bar{D})=\sum_{i=1}^{k}f_{\bar{D}}(H_i)$.
Hence there exists an integer $i$, $1\leq i\leq k$, such that $f_{\bar{D}}(H_i)<c$.
Without loss of generality, assume that $f_{\bar{D}}(H_1)<c$.
Hence $\textrm{cr}_{\bar{D}}(H_1)<c$.
Suppose that $\mathcal{D}=\{D:D$ is a good drawing of $H_1$ with $\textrm{cr}(D)<c\}$.
Let $\mathcal{S}'$ be a subset of $\mathcal{D}$ such that for any $D_1,D_2\in \mathcal{S}'$, $D_1$ is not isomorphic to $D_2$,
and moreover for each $D\in \mathcal{D}$, there exists $D'\in \mathcal{S}'$ such that $D'$ is isomorphic to $D$.
In other words, $\mathcal{S}'$ is the set of all good drawings of $H_1$ up to isomorphism
such that $\textrm{cr}(D')<c$ for all $D'\in \mathcal{S}'$.
Hence $\bar{D}(H_1)\in \mathcal{S}'$, and then $\mathcal{S'}\not=\emptyset$.

Let $1\leq i\leq j\leq k$ in this paragraph.
Let $D_{i,j}$ be a good drawing of $H_{i,j}$.
In $H_{i,j}$, a vertex $u$ is \emph{unsaturated} if $d_{H_{i,j}}(u)<d_{G}(u)$.
If $u,v$ are not on the boundary of a common region in $D_{i,j}$, then $u,v$ are \emph{separated} in $D_{i,j}$.
A $u,v$-\emph{path} is a path with endpoints $u$ and $v$;
the other vertices in a $u,v$-path except for $u,v$ are \emph{internal vertices}.
A $u,v$-path of $D_{i,j}$ is \emph{latent} if both $u$ and $v$ are unsaturated and separated in $D_{i,j}$,
and each internal vertex is not in $V(H_{i,j})$.
A set $\mathcal{P}$ of latent paths of $D_{i,j}$
is \emph{independent} if there is no common edges for any two paths in $\mathcal{P}$.
An independent set $\mathcal{P}$ of latent paths of $D_{i,j}$ is \emph{maximal}
if $\mathcal{P}\cup\{P\}$ is not independent for any latent $u,v$-path $P$ not in $\mathcal{P}$.
Assume that $\emptyset$ is the maximal set of latent paths of $D_{i,j}$ if there are no separated vertices in $D_{i,j}$
or there is no latent $u,v$-path for each pair of separated vertices $u$ and $v$.
Let $D_{i,j}^{+\mathcal{P}}$ be a good drawing of $H_{i,j}\cup \mathcal{P}$ with $D_{i,j}$ as its subdrawing.

Let $\mathcal{S}=\emptyset$.
For $D\in \mathcal{S'}$, let $\mathcal{P}$ be a maximal independent set of latent paths from $D$;
and $\mathcal{S}=\mathcal{S}\cup \{D\}$ if there exists a $D^{+\mathcal{P}}$ such that $f_{D^{+\mathcal{P}}}(H_1)<c$.
For $1\leq i\leq k-1$,
let $\mathcal{D}_i$ be the set of all good drawings of $H_{1,i}$ such that
there exists a $D^{+\mathcal{P}}$ with $f_{D^{+\mathcal{P}}}(H_{1,t})<cj$ for all $D\in \mathcal{D}_i$, $1\leq t\leq i$.

If $\mathcal{D}_1=\emptyset$,
then $l^D_c(H_1)=1$.
Hence $1$=max$\{l^D_c(H_1),...,l^D_c(H_k)\}$.
By Theorem \ref{cor}, we have $\textrm{cr}(G)\geq ck$, a contradiction.
Therefore $\textrm{cr}(G)\geq ck$.

Suppose that $i\geq 2$ and $\mathcal{D}_{i-1}\not=\emptyset$.
Let $\mathcal{S}'=\emptyset$.
For $D\in\mathcal{D}_{i-1}$, we can construct the set $\mathcal{D}'$ of all good drawings of $H_{1,i}$ by adding
$H_i$ to $D$ such that $f_{D'}(H_{1,t})<ct$ for all $D'\in \mathcal{D}'$, $1\leq t\leq i$.
If $i=k$ and $\mathcal{D'}\not=\emptyset$, then $\hat{\mathcal{D}}\not=\emptyset$, and then $\textrm{cr}(G)<ck$.
Otherwise, let $\mathcal{S}'=\mathcal{S}'\cup \mathcal{D'}$.
Let $\mathcal{S}=\emptyset$.
For $D\in \mathcal{S'}$, suppose that $\mathcal{P}$ is a maximal independent set of latent paths of $D$.
Let $\mathcal{S}=\mathcal{S}\cup\{D\}$ if there exists a $D^{+\mathcal{P}}$ such that $f_{D^{+\mathcal{P}}}(H_{1,t})<ct$ for
$1\leq t\leq i$.
If $\mathcal{S}=\emptyset$, then $l^D_c(H_1)=i$.
By a similar argument in the paragraph above,
we have $\textrm{cr}(G)\geq ck$.
Otherwise, let $\mathcal{D}_i=\mathcal{S}$ and let $i=i+1$. Repeat the arguments in this paragraph.
Since $i\leq k$, we can obtain either an integer $t$ such that $\mathcal{D}_t=\emptyset$, or $i=k$ and $\mathcal{D^*}\not=\emptyset$.
By the above analysis, if $\mathcal{D^*}=\emptyset$, then $t$=max$\{l^D_c(H_1),...,l^D_c(H_k)\}$.

To state the exploration more precisely, given a positive number $c$ and a transitive decomposition $\{H_1,...,H_{k}\}$ of a graph $G$,
we introduce an algorithm that outputs $\ell$, which is max$\{l^D_c(H_1),...,l^D_c(H_k)\}$ among all the possible good drawings of $D$, 
if $\textrm{cr}(G)\geq ck$, and a nonempty subset $\mathcal{D^*}$ of $\hat{\mathcal{D}}$, otherwise.
The steps are listed in Algorithm 1.
Algorithm 1 can determine whether a constant $a$ is a lower bound on the crossing number of a given generalized periodic graph.

\begin{table}[!t]
\label{table_time}
\scriptsize
\begin{tabular}{lllllll}
\toprule[2pt]
\textbf{Algorithm 1} Determining-whether-$\textrm{cr}(G)\geq ck$, where $G$ is a generalized periodic graph.\\
\midrule
\textbf{Input:} A positive number $c$, a positive integer $k$, and a generalized periodic graph $G$ with\\

\qquad\quad\ a transitive decomposition $\{H_1,...,H_{k}\}$.\\

\textbf{Output:} $\ell$ or $\mathcal{D^*}$.  \\

\ 1. Construct the set $\mathcal{S'}$ of all good drawings of $H_1$ up to isomorphism such that $\textrm{cr}(D')<c$\\

\quad\, for all $D'\in\mathcal{S'}$.\\

\ 2. \textbf{if} $\mathcal{S'}$ is $\emptyset$ \textbf{then} \\

\ 3.\qquad Let $\ell\gets 1$ and \textbf{return}\\

\ 4. \textbf{end if}\\

\ 5. Let $\mathcal{S}\gets\emptyset$\\

\ 6. \textbf{for all} $D\in \mathcal{S'}$ \textbf{do}\\

\ 7.\qquad Construct a maximal independent set $\mathcal{P}$ of latent paths of $D$.\\

\ 8.\qquad \textbf{if} there exists a $D^{+\mathcal{P}}$ such that $f_{D^{+\mathcal{P}}}(H_1)<c$ \textbf{then}\\

\ 9.\qquad\qquad $\mathcal{S}\gets\mathcal{S}\cup \{D\}$.\\

10.\qquad \textbf{end if}\\

11. \textbf{end for}\\

12. Let $\mathcal{D}_1\gets\mathcal{S}$\\

13. \textbf{if} $\mathcal{D}_1$ is $\emptyset$ \textbf{then}\\

14. \qquad Let $\ell\gets 1$ and \textbf{return}\\

15. \textbf{end if}\\

16. \textbf{for all} $i=2,...,k$ \textbf{do}\\

17. \qquad $\mathcal{S}'\gets\emptyset$.\\

18. \qquad \textbf{for all} $D\in \mathcal{D}_{i-1}$ \textbf{do}\\

19. \qquad\qquad Construct the set $\mathcal{D'}$ of all good drawings of $H_{1,i}$ up to isomorphism by adding\\

 \quad\qquad\qquad\, $H_i$ to $D$ such that $f_{D'}(H_{1,t})<ct$ for all $D'\in\mathcal{D'}$, $1\leq t\leq i$.\\

20. \qquad\qquad \textbf{if} $i$ equals $k$ and $\mathcal{D'}$ is not $\emptyset$ \textbf{then}\\

21. \qquad\qquad\qquad  Let $\mathcal{D^*}\gets\mathcal{D'}$ and \textbf{return}\\

22. \qquad\qquad \textbf{end if}\\

23. \qquad\qquad $\mathcal{S}'\gets\mathcal{S}'\cup\mathcal{D'}$\\

24. \qquad \textbf{end for}\\

25. \qquad $\mathcal{S}\gets\emptyset$.\\

26. \qquad \textbf{for all} $D\in \mathcal{S'}$ \textbf{do}\\

27. \qquad\qquad Construct a maximal independent set $\mathcal{P}$ of latent paths of $D$.\\

28. \qquad\qquad \textbf{if} there exists a $D^{+\mathcal{P}}$ such that $f_{D^{+\mathcal{P}}}(H_{1,t})<ct, 1\leq t\leq i$ \textbf{then}\\

29. \qquad\qquad\qquad $\mathcal{S}=\mathcal{S}\cup \{D\}$.\\

30. \qquad\qquad \textbf{end if}\\

31. \qquad \textbf{end for}\\

32. \qquad Let $\mathcal{D}_i\gets\mathcal{S}$\\

33. \qquad \textbf{if} $\mathcal{D}_i$ is $\emptyset$ \textbf{then}\\

34. \qquad\qquad Let $\ell\gets i$ and \textbf{return}\\

35. \qquad \textbf{end if}\\

36. \textbf{end for}\\
\bottomrule
\end{tabular}
\end{table}
In special cases, by determining a lower bound on the crossing number of a finite generalized periodic graph,
we can determine lower bounds on an infinite family of generalized periodic graphs.
Assume that $G_1,G_2,...$ is an infinite sequence of finite graphs. 
$G_1,G_2,...$ is called \emph{infinite family of generalized periodic graphs on $H$} if
$\{H^i_1,...,H^i_{\phi(i)}\}$ is a transitive decomposition of $G_i,i\geq 1$, where $H$ is a subgraph of $G_1$ 
and $H^i_j$ is isomorphic to $H$ for $1\leq j\leq \phi(i)$, such that $H^i_{1,2}$ is isomorphic to $H^1_{1,2}$ for $i\geq 2$,
and a subdivision of $G_i$ can be obtained by deleting edges from $H^{i+1}_{\phi(i)+1},...,H^{i+1}_{\phi(i+1)}, i\geq 1$.

For example, let $R_i=P(4i+2,2i),i\geq 1$, and let $H$ be the graph $F_0$ defined in Section \ref{P4k2}.
Let $\{F_0,...,F_{\phi(i)-1}\}$ be the transitive decomposition defined in Section \ref{P4k2}.
It is easy to check that $R_1,R_2,...$ is an infinite family of generalized periodic graphs on $F_0$.

We have the following result.
\begin{lemma}\label{a1}
Let $G_1,G_2,...$ be an infinite family of generalized periodic graphs on $H$.
Assume that $\{H^i_1,...,H^i_{\phi(i)}\}$ is a transitive decomposition of $G_i,i\geq 1$, where 
$H^i_j$ is isomorphic to $H$, $1\leq j\leq \phi(i)$, such that $H^i_{1,2}$ is isomorphic to $H^1_{1,2}$, $i\geq 2$,
and a subdivision of $G_i$ can be obtained by deleting edges from $H^{i+1}_{\phi(i)+1},...,H^{i+1}_{\phi(i+1)}, i\geq 1$.
For a given integer $j\geq 1$, Algorithm 1 is used to determine whether $\phi(j)c$ is a lower bound on the crossing number of $G_j$.
Then 

i). If Algorithm 1 returns $\ell$ and $\ell<\phi(j)$,
then $\textrm{cr}(G_i)\geq \phi(i)c$ for $i\geq j$.

ii). If Algorithm 1 returns $\ell$ and $\ell=\phi(j)$,
then a high efficient algorithm Algorithm 2 can determine whether $\textrm{cr}(G_{j+1})\geq \phi(j+1)c$.
\end{lemma}
\begin{proof}
If Algorithm 1 returns $\ell$ and $\ell<\phi(j)$, then $\textrm{cr}(G_j)\geq \phi(j)c$ by Theorem \ref{cor}.
Suppose that $i>j$ and $F$ is a subdivision of $G_j$ obtained by deleting edges from $H^i_{\phi(j)+1},...,H^i_{\phi(i)}$.
Since $\{H^i_1,...,H^i_{\phi(i)}\}$ and $\{H^j_1,...,H^j_{\phi(j)}\}$ are transitive decompositions of $G_i$ and $G_j$, respectively,
both $H^i_k$ and $H^j_t$ are isomorphic to $H$ for $1\leq k\leq \phi(i)$ and $1\leq t\leq \phi(j)$, $H^i_{1,2}$ is isomorphic to $H^j_{1,2}$,
each step of the progress that Algorithm 1 determining whether $\textrm{cr}(G_j)\geq \phi(j)c$ is always valid in $F$ by the construction of $F$.
Hence max$\{l^D_c(H^i_1),...,l^D_c(H^i_{\phi(i)})\}=$max$\{l^D_c(H^j_1),...,$ $l^D_c(H^j_{\phi(j)})\}<\phi(j)<\phi(i)$, and then $\textrm{cr}(G_i)\geq \phi(i)c$ by Theorem \ref{cor}.

Suppose that Algorithm 1 returns $\ell$ and $\ell=\phi(j)$.
Then $\textrm{cr}(G_j)\geq \phi(j)c$ by Theorem \ref{cor}.
Let $\mathcal{E}$ be the set of set $E$ of edges such that a subdivision of $G_j$ can be obtained by deleting $E$ from $G_{j+1}$.
Then $\mathcal{E}$ is not empty since $E(G_{j+1})-E(F)$ is in $\mathcal{E}$.
Moreover, $|\mathcal{E}|\geq \phi(j+1)$ since $\{H^{j+1}_1,...,H^{j+1}_{\phi(j+1)}\}$ is a transitive decomposition of $G_{j+1}$.
Assume that $|\mathcal{E}|=\psi(j+1)$ and $\mathcal{E}=\{E_1,...,E_{\psi(j+1)}\}$.

Suppose that $\textrm{cr}(G_{j+1})<\phi(j+1)c$ and
$D$ is a good drawing of $G_{j+1}$ with $\textrm{cr}(D)<\phi(j+1)c$.
Let $d=\phi(j+1)-\phi(j)$.
If there exists an integer $i$, $1\leq i\leq \psi(j+1)$, such that $\textrm{cr}_{D}(E_i)+\textrm{cr}_{D}(E_i,G_{j+1}\setminus E_i)\geq cd$, 
then a good drawing $D'$ of a subdivision of $G_j$ with $\textrm{cr}(D')<\phi(j)c$ can obtained by deleting $E_i$ from $D$,
which contradicts $\textrm{cr}(G_j)\geq \phi(j)c$.
Hence $\textrm{cr}_{D}(E_i)+\textrm{cr}_{D}(E_i,G_{j+1}\setminus E_i)<cd$, $1\leq i\leq \psi(j+1)$.
Therefore, a high efficient algorithm can be obtained by modifying several steps in Algorithm 1 to determine whether $\textrm{cr}(G_{j+1})\geq\phi(j+1)c$ with the assumption that $\textrm{cr}(G_j)\geq\phi(j)c$.
The modification is listed in Algorithm 2.
Algorithm 2 outputs $\ell$ if $\textrm{cr}(G_{j+1})\geq\phi(j+1)c$,
and a nonempty set $\mathcal{D^*}$ of good drawings $D$ of $G_{j+1}$ with $\textrm{cr}(D)<\phi(j+1)c$, otherwise.
\end{proof}
\begin{remark}\label{rem2}
The proof of Theorem \ref{thm1} is an instance of Algorithm 2.
\end{remark}
\begin{table}[!t]
\label{table_time}
\scriptsize
\begin{tabular}{lllllll}
\toprule[2pt]
\textbf{Algorithm 2} Determining-whether-$\textrm{cr}(G_{j+1})\geq \phi(j+1)c$-with-$\textrm{cr}(G_j)\geq \phi(j)c$\\
\midrule
\textbf{Input:} A positive number $c$, positive integers $k:=\phi(j+1)$, $\psi(j+1)$ and $d:=\phi(j+1)-\phi(j)$, \\

\qquad\quad\ a generalized periodic graph $G_{j+1}$ with a transitive decomposition $\{H^{j+1}_1,...,H^{j+1}_{\phi(j+1)}\}$,\\

\qquad\quad\ and $\{E_1,...,E_{\psi(j+1)}\}$ such that $G_{j+1}\setminus E_i$ is a subdivision of $G_j,1\leq i\leq \psi(j+1)$.\\

\textbf{Output:} $\ell$ or $\mathcal{D^*}$  \\

The same steps (with symbol $H$ replaced by $H^{j+1}$) as those of Algorithm 1 except for 4 lines;\\

those which line $i$ in Algorithm 1 is replaced by step $i'$ as follows, where $i\in\{1,8,19,28\}$.\\

\ 1'. Construct the set $\mathcal{S'}$ of all good drawings of $H^{j+1}_1$ up to isomorphism such that $\textrm{cr}(D')<c$\\

\qquad and $\textrm{cr}_{D'}(H^{j+1}_1\cap E_i)+\textrm{cr}_{D'}(H^{j+1}_1\cap E_i,G_{j+1}\setminus(H^{j+1}_1\cap E_i))<cd$ for all $D'\in\mathcal{S'}$ and for \\

\qquad all $H^{j+1}_1\cap E_i\neq\emptyset, 1\leq i\leq \psi(j+1)$.\\

\ 8'. \textbf{if} there exists a $D^{+\mathcal{P}}$ such that $f_{D^{+\mathcal{P}}}(H^{j+1}_1)<c$ and $\textrm{cr}_{D^{+\mathcal{P}}}(H^{j+1}_1\cap E_i)+\textrm{cr}_{D^{+\mathcal{P}}}(H^{j+1}_1\cap$\\

\qquad $E_i,G_{j+1}\setminus(H^{j+1}_1\cap E_i))<cd$ for all $H^{j+1}_1\cap E_i\neq\emptyset, 1\leq i\leq \psi(j+1)$ \textbf{then}\\

19'.  Construct the set $\mathcal{D'}$ of all good drawings of $H^{j+1}_{1,i}$ up to isomorphism by adding $H^{j+1}_i$ to\\

\qquad $D$ such that $f_{D'}(H^{j+1}_{1,t})<ct$, $1\leq t\leq i$, and $\textrm{cr}_{D'}(H^{j+1}_{1,i}\cap E_t)+\textrm{cr}_{D'}(H^{j+1}_{1,i}\cap E_t,H^{j+1}_{1,i}\setminus$\\

\qquad $(H^{j+1}_{1,i}\cap E_t))<cd$ for all $H^{j+1}_{1,i}\cap E_t\neq\emptyset, 1\leq t\leq \psi(j+1)$ and for all $D'\in\mathcal{D'}$.\\

28'. \textbf{if} there exists a $D^{+\mathcal{P}}$ such that $f_{D^{+\mathcal{P}}}(H^{j+1}_{1,t})<ct, 1\leq t\leq i$, and $\textrm{cr}_{D'}(H^{j+1}_{1,i}\cap E_t)+$\\

\qquad $\textrm{cr}_{D'}(H^{j+1}_{1,i}\cap E_t,H^{j+1}_{1,i}\setminus(H^{j+1}_{1,i}\cap E_t))<cd$ for all $H^{j+1}_{1,i}\cap E_t\neq\emptyset, 1\leq t\leq \psi(j+1)$ \textbf{then}\\
\bottomrule
\end{tabular}
\end{table}

\section{Acknowledgement.}

The work was supported by the National Natural Science Foundation of China (No.\,61401186, No.\,11901268).

\section{Reference}







\end{document}